\documentclass[12pt,
reqno]{amsart}

\usepackage[cp1251]{inputenc}
\usepackage{amsmath,amsxtra,amssymb,eufrak}
\usepackage[all]{xy}
\usepackage[active]{srcltx} 
\usepackage{mathrsfs}
\usepackage[russian]{babel}

\theoremstyle{plain}
\newtheorem{thm}{Теорема}[section]
\newtheorem{lm}[thm]{Лемма}
\newtheorem{co}[thm]{Следствие}
\newtheorem{pr}[thm]{Предложение}

\theoremstyle{definition}

\newtheorem{df}[thm]{Определение}
\newtheorem{exm}[thm]{Пример}
\newtheorem{rem}[thm]{Замечание}
\newtheorem{que}[thm]{Вопрос}

\numberwithin{equation}{section}

\newcommand*{\Ker}{\mathop{\mathrm{Ker}}}
\newcommand*{\Image}{\mathop{\mathrm{Im}}}
\newcommand*{\Smash}{\mathop{\#}}
\newcommand*{\Psmash}{\mathbin{\widehat{\#}}}
\newcommand*{\ptn}{\mathbin{\widehat{\otimes}}}
\newcommand*{\ptens}[1]{\mathop{\widehat\otimes}_{#1}}
\newcommand*{\lmod}{\mbox{-}\!\mathop{\mathsf{mod}}}
\newcommand*{\rmod}{\mathop{\mathsf{mod}}\!\mbox{-}}
\newcommand*{\bimod}{\mbox{-}\!\mathop{\mathsf{mod}}\!\mbox{-}}

\title{Алгебры аналитических функционалов и гомологические эпиморфизмы}

\author{О. Ю. Аристов}
\email{aristovoyu@inbox.ru}
\address{Обнинск, Калужская Область}

\begin{document}
\begin{abstract}
В статье автора [arXiv: 2404.19433] было установлено, что оболочка Аренса--Майкла разрешимой алгебры Ли является гомологическим эпиморфизмом. Здесь показано, что для алгебр аналитических функционалов на связной комплексной группе Ли аналогичное утверждение выполнено без предположения о разрешимости, и, более того, гомоморфизмы пополнения более общего вида также являются гомологическими эпиморфизмами. В частности, это выполнено для оболочки относительно класса банаховых PI-алгебр. 
\end{abstract}

 \maketitle
 \markright{Аналитические функционалы и гомологические эпиморфизмы}

\section{Введение}

Оба основных понятия, обсуждаемые в этой статье, ---  гомологические эпиморфизмы и оболочки Аренса--Майкла  --- были введены Джозефом Тейлором в его работах начала 1970-х (хотя и под другими названиями). Им же в неявной форме был поставлен вопрос: \emph{при каких условиях для данной конечномерной комплексной алгебры Ли~$\mathfrak{g}$ оболочка Аренса--Майкла  её универсальной обертывающей алгебры $U(\mathfrak{g})$ является гомологическим эпиморфизмом?} (Явно сформулирован в \cite[\S\,9, Problem~1]{Pir_stbflat}). В результате ряда усилий, приложенных Доси, Пирковским и автором, окончательный ответ был получен в~\cite{Ar_smhe}:  указанная \emph{оболочка является гомологическим эпиморфизмом тогда и только тогда, когда $\mathfrak{g}$ разрешима}.

Основной целью этой работы является доказательство того, что ограничение на разрешимость отсутствует
для алгебры аналитических функционалов ${\mathscr A}(G)$ на связной комплексной группе Ли~$G$. А именно, мы покажем, что \emph{оболочка Аренса--Майкла  ${\mathscr A}(G)\to \widehat{\mathscr A}(G)$ и оболочка относительно класса банаховых PI-алгебр ${\mathscr A}(G)\to\widehat{\mathscr A}(G)^{\mathrm{PI}}$ являются гомологическими эпиморфизмами} (теорема~\ref{C4main}). 

Эта статья является завершением серии работ автора, посвящённых темам, тесно связанным с рассматриваемой задачей. Содержащиеся в них утверждения, помимо того, что необходимы для доказательства теоремы~\ref{C4mainexdi} (общего утверждения, из которого выводятся основные результаты этой статьи), также представляют и самостоятельный интерес. Разложение алгебры аналитических функционалов и её пополнений (в том числе $\widehat{\mathscr A}(G)$ и $\widehat{\mathscr A}(G)^{\mathrm{PI}}$) в итерированные аналитические смэш-произведения получены в \cite{Ar_smash}, а обсуждение гомологических эпиморфизмов для универсальных обертывающих алгебр содержится в \cite{Ar_smhe}. Первая из этих статей опирается на результаты  о разложении функций длины из \cite{ArLfd} (которая, в свою очередь, обобщает \cite{ArAnF}) и \cite{ArPiLie}, где изучен вопрос о том, когда пополнение универсальной обертывающей алгебры является банаховой PI-алгеброй. Кроме того, в \cite{Ar_lin} было уточнено строение линеаризатора связной группы Ли, что также необходимо для доказательства  теоремы~\ref{C4mainexdi}.

\subsection*{Гомологические эпиморфизмы}
В некоторых вопросах алгебры и анализа требование того, что гомоморфизм алгебр является плоским, оказывается слишком сильным и должно быть смягчено. Соответствующее более слабое понятие в контексте топологических алгебр формулируется так. Непрерывный гомоморфизм полных локально выпуклых алгебр с совместно непрерывным умножением (мы называем их $\ptn$-алгебрами) называется \emph{гомологическим эпиморфизмом}, если  производный функтор от функтора отступления между соответствующими категориями модулей является вполне унивалентным. В этой статье язык производных категорий оказался не нужен, так как достаточно использовать эквивалентное определение в классических терминах резольвент, см. \S\,\ref{sec:prel}. Само понятие гомологического эпиморфизма восходит к Дж. Тейлору \cite{T2} и было неоднократно переоткрыто в разных контекстах и под разными именами, см. обсуждение и историю терминологии в \cite[Remark 3.16]{AP}.

Важным для приложений свойством является то, что гомологические эпиморфизмы сохраняют непрерывные когомологии с коэффициентами в $\ptn$-модулях \cite[Proposition 1.7(a)]{T2}. Например, если естественное отображение между групповыми алгебрами различного типа  является гомологическим эпиморфизмом, то соответствующие когомологии совпадают, см., например, частные случаи в \cite{Ji92,Me06}.

Естественным приложением гомологических эпиморфизмов является оценка гомологических размерностей, см. \cite{Ar_dim}. Также в качестве приложения можно привести препринт Майера~\cite{Me04}, где показано, что, используя их, можно не только упростить  вычисление циклических когомологий гладкого некоммутативного тора, предпринятое Конном в \cite[\S\,6]{Co85}, но и сделать его более понятным. Кроме того, слабая версия понятия гомологического эпиморфизма применена Пирковским для характеризации открытых вложений многообразий Штейна~\cite{Pi99}, а в дальнейшем этот результат был распространён и на общие пространства Штейна \cite{AP,BBK18}.

\subsection*{Связь со спектральной теорией}
Обращаясь к исходному вопросу о свойствах оболочек Аренса--Майкла, приведём мотивировку --- откуда возникла потребность изучать их с гомологической точки зрения. Проблематика появилась в работах Тейлора начала 70-ых \cite{T70,T1,T2,Ta73}, посвящённых спектральной теории нескольких коммутирующих и некоммутирующих операторов в банаховых пространствах. Тейлор предложил гомологическое определение спектра представления (модуля), основанное на производном функторе ${\mathop{\mathrm{Tor}}\nolimits}$. Пусть $A$ --- некоторая алгебра над~$\mathbb{C}$, например, конечнопорождённая (что соответствует случаю нескольких операторов, удовлетворяющих заданным соотношениям), а вообще говоря, локально выпуклая. Спектр левого локально выпуклого модуля $N$  определяется как подмножество каким-либо образом выбранного структурного пространства (состоящего из классов эквивалентности правых модулей~$M$) заданное условием ${\mathop{\mathrm{Tor}}\nolimits}_n^{A}(M,N)\neq 0$ для некоторого $n\in\mathbb{Z}_+$ (здесь мы следуем интерпретации идей Тейлора из введения к \cite{Pir_stbflat}, см. там же подробное обсуждение). Отметим, что для некоммутативных алгебр желательно ещё установить соответствие между правыми и левыми модулями.

Уже несколько десятилетий развивается спектральная теория представлений разрешимых алгебр Ли в банаховых пространствах (основанная на одномерных представлениях), инициированная Файнштейном \cite{Fa93} (для нильпотентных), а также Боассо и Ларотонда \cite{BL93}, см. книгу \cite{BS01}. Она использует более простое определение спектра через комплекс Кошуля, идея которого тоже восходит к Тейлору. Из более современных работ следует указать серию статей Доси, см., например, \cite{Do05,Do09,Do10C,Do10A}.

Основным мотивом Тейлора при введении гомологического понятия спектра было, по-видимому, желание построить функциональное исчисление не только для коммутирующих операторов (в чём он преуспел), но и для некоммутирующих, где он только заложил основы теории. Надо признать, что, в отличие от коммутативной, успехи некоммутативной спектральной теории в построении функциональных исчислений не слишком велики. Наиболее сильным результатом является частное  (но при этом довольно трудное) утверждение о функциональном исчислении для нильпотентных алгебр Ли \cite{Do10A}. Поведение же спектра,  определённого через Tor на множестве одномерных представлений, имеет патологии даже в простейших примерах разрешимых алгебр Ли, не являющихся нильпотентными, см., например, \cite{Bi20}.

Идея использовать при определении спектра элементы некоторой большей алгебры восходит к классическим образцам. Ещё во введениях к статьям \cite{T2,Ta73} Тейлор обратил внимание на важность тех гомоморфизмов в некоторую алгебру ''некоммутативных аналитических функций'',  которые  являются гомологическими в современной терминологии. Гомоморфизм $A\to B$ индуцирует последовательность линейных отображений
$$
{\mathop{\mathrm{Tor}}\nolimits}_n^{A}(M,N) \to {\mathop{\mathrm{Tor}}\nolimits}_n^{B}(M,N)\qquad (n\in\mathbb{Z}_+)
$$
для соответственно правого и левого $B$-$\ptn$-модуля $M$ и $N$.
Если эти отображения биективны для всех $M$, $N$ и $n$ (т.е. гомоморфизм сохраняет кручение, см. \cite[Добавление]{He81}), то, очевидно, данное выше определение спектра не зависит от того, какая алгебра использована, $A$ или $B$.  Для биективности достаточно, чтобы $A\to B$ был гомологическим эпиморфизмом \cite[Proposition 1.4]{T2}.

Что касается выбора алгебры $B$, то если речь идёт от представлениях в банаховых пространствах, её можно считать банаховой или по крайней мере, аппроксимируемой банаховыми (т.е. алгеброй Аренса--Майкла). Естественным кандидатом оказывается универсальная алгебра Аренса--Майкла, ассоциированная с $A$, --- её оболочка Аренса--Майкла $\widehat A$ (в силу того, что структура банахова модуля автоматически продолжается с $A$  на $\widehat A$). Напомним, что оболочка Аренса--Майкла топологической алгебры  может быть определена как пополнение относительно семейства всех непрерывных субмультипликативных преднорм. Она обладает универсальным свойством в классе всех банаховых алгебр. Ограничивая класс алгебр, мы получаем также оболочку в классе банаховых PI-алгебр (обозначается через $\widehat A^{\mathrm{PI}}$). Это понятие  введено в \cite{Ar_smash} при исследовании случая $A={\mathscr A}(G)$. (Напомним, что ${\mathscr A}(G)$ --- это свёрточная алгебра аналитических функционалов на комплексной группе Ли $G$,  подробности см. в \S\,\ref{sec:prel}.)

Таким образом, вопрос о том, когда оболочка Аренса--Майкла является гомологическим эпиморфизмом, тесно связан со спектральной теорией. В частности, рассматривая вопрос о том, будут ли эпиморфимы ${\mathscr A}(G)\to \widehat{\mathscr A}(G)$ и ${\mathscr A}(G)\to \widehat{\mathscr A}(G)^{\mathrm{PI}}$ являться гомологическими, мы сравниваем спектральные теории голоморфных представлений $G$ в локально выпуклых и в банаховых пространствах  (во втором случае --- с дополнительным условием наличия полиномиального тождества, общего для всех элементов образа)\footnote{Отметим, что алгебра $\widehat{\mathscr A}(G)$  имеет и другие приложения. Например, переход от ${\mathscr A}(G)$ к $\widehat{\mathscr A}(G)$ является существенным шагом в конструкции голоморфной двойственности, предложенной Акбаровым в \cite{Ak08}, см. также \cite{ArHR}.}. Однако настоящая статья не содержит результатов по спектральной теории --- для нас она является лишь источником исходного вопроса.

\subsection*{Алгебра аналитических функционалов и её пополнения}

Если $G$ односвязна, а $\mathfrak{g}$ --- её алгебра Ли, то $\widehat{\mathscr A}(G)\cong \widehat U(\mathfrak{g}) $ и $\widehat{\mathscr A}(G)^{\mathrm{PI}}\cong \widehat U(\mathfrak{g})^{\mathrm{PI}}$. Поэтому изучать $\widehat{\mathscr A}(G)$ и $\widehat{\mathscr A}(G)^{\mathrm{PI}}$ как  обобщения соответственно $\widehat U(\mathfrak{g})$ и $\widehat U(\mathfrak{g})^{\mathrm{PI}}$ кажется довольно естественным.

В \cite{Ar_smash} показано, что $\widehat{\mathscr A}(G)$ и $\widehat{\mathscr A}(G)^{\mathrm{PI}}$  являются крайними членами семейства пополнений алгебры ${\mathscr A}_{\omega^\infty}(G)$, индексированного субмультипликативными весами $\omega$, обладающими некоторым свойством максимальности. Как шаг к основному результату,  мы докажем общее утверждение, включающее $\widehat{\mathscr A}(G)$ и $\widehat{\mathscr A}(G)^{\mathrm{PI}}$  как частные случаи: если $G$ --- связная линейная комплексная группа Ли, то ${\mathscr A}(G)\to {\mathscr A}_{\omega^\infty}(G)$ является гомологическим эпиморфизмом (см. еще более общую формулировку в теореме~\ref{C4mainexdi}).

\subsection*{Смэш-произведения}
Для того чтобы доказать, что ${\mathscr A}(G)\to {\mathscr A}_{\omega^\infty}(G)$ --- гомологический эпиморфизм, мы будем использовать согласованные разложения обеих алгебр в итерированное аналитическое смэш-произведение, полученные в~\cite{Ar_smash}.
Понятие аналитического смэш-произведения ввёл Пирковский в~\cite{Pi4} именно с целью изучения гомологических свойств оболочек Аренса--Майкла.

Результаты автора из \cite{Ar_smash} вкратце сводятся к следующему. Для связной линейной группы Ли $G$ можно построить разложение в итерированное полупрямое произведение
$$
G\cong ((\cdots (F_1 \rtimes F_2)\rtimes\cdots)
\rtimes F_n),
$$
таким образом, чтобы не только ${\mathscr A}(G)$  разлагалось в итерированное смэш-произведение
$$
((\cdots ({\mathscr A}(F_1) \Psmash {\mathscr A}(F_2))\Psmash\cdots)
\Psmash {\mathscr A}(F_n))
$$
(что не так уж трудно показать), но также и ${\mathscr A}_{{\omega}^\infty}(G)$:
$$
{\mathscr A}_{{\omega}^\infty}(G)\cong((\cdots ({\mathscr A}_{\omega_1^\infty}(F_1) \mathop{\widehat{\#}} {\mathscr A}_{\omega_2^\infty}(F_2))\mathop{\widehat{\#}}\cdots)
\mathop{\widehat{\#}} {\mathscr A}_{\omega_{k+1}^\infty}(F_{k+1}));
$$
здесь ${\mathscr A}_{\omega_i^\infty}(F_i)$ --- специальным образом выбранное пополнение ${\mathscr A}(F_i)$.
Существенно, что эпиморфизмы ${\mathscr A}(F_i)\to {\mathscr A}_{\omega_i^\infty}(F_i)$ являются гомологическими и  связаны между собой дополнительными гомологическими соотношениями, см. \S\,\ref{s:proofmaingen}.

\subsection*{Редукция к разрешимому случаю}

Общая связная комплексная группа Ли разлагается в композиционный ряд
$$
1\lhd  \mathrm{M}\lhd \Lambda \lhd B \lhd G,
$$ где $\mathrm{M}$ --- подгруппа, на которой все голоморфные функции постоянны, $\Lambda$ --- линеаризатор,
$B/\Lambda$ односвязна и разрешима, а $G/B$ линейно комплексно редуктивна.

В силу определения ${\mathscr A}(G)$ мы можем игнорировать $\mathrm{M}$. Группы $\Lambda/ \mathrm{M}$ и $G/B$ линейно комплексно редуктивны, что влечёт гомологическую тривиальность ${\mathscr A}(\Lambda)$ и ${\mathscr A}(G/B)$. В частности, это означает, что мы можем игнорировать  $\Lambda$ при построении резольвенты. Хотя $G/B$ игнорировать не получится, её нетривиальность по крайней мере не отражается на гомологических характеристиках. Таким образом, вся нужная информация о  ${\mathscr A}(G)$ и её пополнениях может быть в итоге выведена из свойств~$B/\Lambda$. 
Поскольку $B/\Lambda$ односвязна и разрешима, задача сводится к изучению пополнений универсальной обертывающей алгебры разрешимой группы Ли. Это в основном проделано в работе автора~\cite{Ar_smhe}, см. там же историю этого вопроса.

\subsection*{Открытые вопросы}
Следующие открытые вопросы представляются интересными.

Понятие локализующего веса на комплексной группе Ли введено ниже в определении~\ref{deflocwe}.

\begin{que}\label{locwge}
При каких общих условиях  вес является локализующим?
\end{que}

Пусть теперь $G$ --- связная действительная группа Ли. Обозначим через ${\mathscr E}'(G)$ алгебру распределений с компактным носителем (т.е. непрерывных линейных функционалов на пространстве всех гладких функций) со сверткой.

\begin{que}\label{EGprHE}
При каких условиях ${\mathscr E}'(G)\to\widehat{\mathscr E}'(G)$  и ${\mathscr E}'(G)\to\widehat{\mathscr E}'(G)^{\mathrm{PI}}$ являются гомологическими эпиморфизмами?
\end{que}

\subsection*{Благодарности}
Автор хотел бы выразить признательность рецензенту за два высококачественных отзыва, содержащих многочисленные замечания, учёт которых значительно улучшил изложение (в частности были упрощены доказательства предложений~\ref{norbal2}, \ref{hotrqufr}  и~\ref{qLaloca}),  а также за  ссылку с доказательством леммы~\ref{topxyprxpr}.

Большая часть этой работы была выполнена во время пребывания автора в должности исследователя в Институте перспективных исследований в области математики Харбинского Политехнического Университета. Автор благодарен указанному учреждению за поддержку. Кроме того, часть работы выполнена во время визита автора в НИУ ВШЭ (г.\,Москва) зимой 2025\,г. Автор выражает признательность этому университету и особенно А.\,Ю. Пирковскому за гостеприимство. Автор также благодарен  А.\,B.~Домрину за полезные консультации и участникам семинара ``Алгебры в Анализе'' (мехмат МГУ), где предварительная версия этой работы была рассказана в 2018 году. 

\tableofcontents

\section{Предварительные определения и обозначения}
\label{sec:prel}

\subsection*{Аналитические функционалы}

Для комплексного многообразия $M$ мы обозначаем через $\mathcal{O}(M)$
локально выпуклое пространство  всех голоморфных функций на~$M$ (с топологией равномерной сходимости на компактных подмножествах), а через  ${\mathscr A}(M)$ --- сильное двойственное пространство $\mathcal{O}(M)'$  (множество непрерывных линейных функционалов с топологией равномерной сходимости на ограниченных подмножествах). Если $N$ и $M$ --- комплексные многообразия, то поскольку $\mathcal{O}(N)$ и
$\mathcal{O}(M)$ являются ядерными пространствами Фреше, имеют место следующие топологические изоморфизмы  полных локально выпуклых пространств (они понадобятся ниже)
\begin{multline}
\label{LSCtnpr}
 {\mathscr A}(N\times M)=\mathcal{O}(N\times M)'\cong\\
\cong(\mathcal{O}(N)\ptn\mathcal{O}(M))'\cong  \mathcal{O}(N)'\ptn\mathcal{O}(M)'= {\mathscr A}(N)\ptn {\mathscr A}(M).
\end{multline}
(Здесь $\ptn$ обозначает полное проективное тензорное произведение локально выпуклых пространств.)

Пусть теперь $G$  --- комплексная группа Ли. Тогда $\mathcal{O}(G)$ является не только
$\ptn$-алгеброй  (относительно поточечного умножения), но и $\ptn$-алгеброй Хопфа.
(Напомним, что \emph{$\ptn$-алгеброй Хопфа} называется алгебра Хопфа
в симметрической моноидальной категории полных локально выпуклых пространств с бифунктором $(-)\ptn(-)$. Аналогично определяются $\ptn$-коалгебры и $\ptn$-биалгебры.)

Формулы
\begin{equation*}
 \Delta_\mathcal{O} (f)(g, h) = f(gh),\quad \varepsilon_\mathcal{O} (f) = f(1),\quad  (S_\mathcal{O} f)(g) = f(g^{-1})
\end{equation*}
задают на $\mathcal{O}(G)$ каноническую структуру $\ptn$-алгебры Хопфа.
Тогда ${\mathscr A}(G)$ можно снабдить двойственной структурой $\ptn$-алгебры Хопфа.
Умножение (свёртка) на ${\mathscr A}(G)$  задаётся формулой
$$
\langle \alpha\beta, f\rangle\!:= (\alpha\otimes \beta)\Delta_\mathcal{O} (f) \qquad(\alpha\in{\mathscr A}(G),\,f\in\mathcal{O}(G)),
$$
а единицей является дельта-функция $\delta_e$ в единице группы. Остальные операции определяются формулами
$$
\Delta_{\mathscr A}(\alpha)(f_1\otimes f_2)=\langle \alpha, f_1f_2\rangle,\quad\varepsilon_{\mathscr A}(\alpha)=\langle \alpha, 1\rangle,\quad \langle S_{\mathscr A}(\alpha), f\rangle = \langle \alpha, S_\mathcal{O}(f)\rangle
$$
(здесь $\alpha\in{\mathscr A}(G)$ и $f,f_1,f_2\in\mathcal{O}(G)$).

\subsection*{Cубмультипликативные веса}
Напомним, что строго положительная локально ограниченная функция $\omega\!: G \to\mathbb{R}$  на локально компактной группе $G$ называется \emph{субмультипликативным весом}, если $\omega(gh)\le \omega(g)\omega(h)$ для всех $g, h \in G$.

Мы говорим, что вес $\omega_2$ \emph{мажорирует} вес $\omega_1$ (пишем
$\omega_1\preceq \omega_2$), если найдутся $C>0$ и $\gamma>0$, такие что
$$
\omega_1(g)\le C\,\omega_2(g)^\gamma \quad\text{для всех $g\in G$.}
$$
Если $\omega_1\preceq \omega_2$ и $\omega_2\preceq \omega_1$, то мы называем веса эквивалентными (пишем $\omega_1\simeq  \omega_2$).

Субмультипликативный вес $\omega$ на локально компактной группе $G$ называется \emph{асимптотически симметричным}, если $\omega$ мажорирует функцию $g\mapsto\omega(g^{-1})$ и $\omega(g)\ge 1$ для всех $g\in G$ \cite[определение 3.5]{Ar_smash}.

Хорошо известно, что на компактно порождённой локально компактной группе существует субмультипликативный вес, который мажорирует все субмультипликативные веса на ней; в частности, это верно для связных групп. Мы будем назвать его \emph{максимальным весом}.

Пусть $\omega$ --- субмультипликативный вес на локально компактной группе $G$, а $\omega_H$ --- максимальный вес на её замкнутой подгруппе~$H$.  Следуя \cite{Ar_smash}, будем говорить, что $\omega$ имеет \emph{экспоненциальное искривление
} на $H$, если найдутся $C>0$ и $\gamma>0$, такие что $1+\log \omega_H$ мажорирует $\omega$ на $H$.

Пусть $\mathfrak{g}$ --- алгебра Ли связной линейной группы  Ли $G$. Напомним, что подгруппа, порождённая $\exp\mathfrak{h}$ для некоторой подалгебры Ли $\mathfrak{h}$ в $\mathfrak{g}$ называется \emph{интегральной}. Обозначим через $N$ \emph{нильпотентный  радикал}~$G$ --- интегральную подгруппу, алгебра Ли которой есть нильпотентный  радикал алгебры Ли  $\mathfrak{g}$, т.е. пересечение ядер всех конечномерных неприводимых представлений~$\mathfrak{g}$ (см. \cite[глава~I, \S\,5.3, с.\,58, определение~3]{Bou} или \cite[\S\,1.7.2]{Dix}).
Обозначим также через  $E$ \emph{экспоненциальный  радикал}~$G$. Мы опускаем его определение, сам термин  предложен в~\cite{Os02}, а в нужной нам общности оно дано в~\cite{Co08}; см. также \cite[\S\,3]{ArAnF}.

Если $N'$ ---  нормальная интегральная подгруппа  такая, что $E\subset N'\subset N$, то существует субмультипликативный вес, который мажорирует все субмультипликативные веса, имеющие экспоненциальное искривление на $N'$, см. \cite[Theorem 2.3]{ArLfd}  и \cite[Theorem 4.3]{Ar_smash}. 
Мы называем его \emph{максимальным весом с экспоненциальным искривлением на~$N'$},

\subsection*{Пополнения алгебры аналитических функционалов}
Пусть теперь $G$ --- комплексная группа Ли, а $\omega$ --- субмультипликативный вес на~$G$. Обозначим замыкание в~${\mathscr A}(G)$ абсолютно выпуклой оболочки подмножества
\begin{equation}\label{Vup}
\{\omega(g)^{-1}\delta_g:\,x\in G\}
\end{equation}
через $V_\omega$ (здесь $\delta_g$ обозначает дельта-функцию в точке $g$).  Функционал Минковского на ${{\mathscr A}}(G)$, ассоциированный с $V_\omega$, является непрерывной полунормой на ${{\mathscr A}}(G)$. Обозначим его через $\|\cdot\|_{\omega}$,  а через ${{\mathscr A}}_{\omega^\infty}(G)$ --- пополнение ${{\mathscr A}}(G)$ относительно последовательности преднорм 
$$
(\|\cdot\|_{\omega^n};\,n\in\mathbb{N}),\qquad \text{где $\omega^n(g)\!:=\omega(g)^n$.}
$$
Если $\omega\ge 1$, то умножение на ${\mathscr A}(G)$ продолжается до умножения на ${\mathscr A}_{\omega^\infty}(G)$, превращающего её в алгебру Фреше--Аренса--Майкла \cite[предложение~3.3]{ArPiLie}, 
см. также~\cite{Ak08}, где эта конструкция была впервые введена.

\subsection*{Оболочки}
Пусть $\mathcal{C}$ --- некоторый класс банаховых алгебр. Мы говорим, что топологическая алгебра \emph{локально содержится в классе $\mathcal{C}$}, если она является проективным пределом алгебр из $\mathcal{C}$ в категории топологических алгебр.

\begin{df}\label{Bcalssun}
Пусть $A$ --- ассоциативная топологическая алгебра, а $\mathcal{C}$ --- класс банаховых алгебр, стабильный относительно перехода к замкнутым подалгебрам и конечным произведениям. Мы говорим, что пара $(\widehat A^{\mathcal{C}},\iota)$, где $\widehat A^{\mathcal{C}}$ -локально содержится в $\mathcal{C}$, а $\iota$ непрерывный гомоморфизм $A\to\widehat A^{\mathcal{C}}$, является \emph{оболочкой $A$ в классе $\mathcal{C}$}, если выполнено следующее универсальное свойство.
Для всякой $B\in\mathcal{C}$ и всякого непрерывного гомоморфизма
$\varphi\!: A\to B$ существует единственный непрерывный гомоморфизм
$\widehat\varphi\!:\widehat A^{\mathcal{C}}\to B$, такой что диаграмма
\begin{equation}\label{AMen}
  \xymatrix{
A \ar[r]^{\iota}\ar[rd]_{\varphi}&\widehat A^\mathcal{C}\ar@{-->}[d]^{\widehat\varphi}\\
 &B\\
 }
\end{equation}
коммутативна.
\end{df}
Нетрудно показать, что  существование оболочки  следует из требования стабильности $\mathcal{C}$ относительно перехода к замкнутым подалгебрам и конечным произведениям, см. подробности в~\cite{Ar_envPG}. Как легко видеть, оболочка единственна с точностью до изоморфизма. Из определения сразу следует, что оболочка является функтором, т.е. всякий непрерывный гомоморфизм $\theta\!:A\to B$ топологических алгебр порождает непрерывный гомоморфизм $\widehat{\theta}^\mathcal{C}\!:\widehat A^\mathcal{C}\to \widehat B^\mathcal{C}$.

Если $\mathcal{C}$ --- класс всех банаховых алгебр, то получаем хорошо известную \emph{оболочку Аренса--Майкла}, которая обозначается через $\widehat A$. Если $\mathcal{C}$ --- класс  банаховых PI-алгебр, то соответствующую оболочку обозначаем через $\widehat A^{\mathrm{PI}}$ (она впервые определена в \cite{Ar_smash}). (Напомним, что ассоциативная алгебра $A$  называется PI-\emph{алгеброй}, если существуют $n\in\mathbb{N}$ и ненулевой элемент $p$ свободной алгебры с $n$ образующими, такие что  $p(a_1,\ldots,a_n)=0$ для всех $a_1,\ldots,a_n$, см., например, \cite{KKR16}.)

Нас интересуют оболочки алгебры ${\mathscr A}(G)$, где $G$ --- комплексная группа Ли. В связном случае $\widehat{\mathscr A}(G)$ и $\widehat{\mathscr A}(G)^{\mathrm{PI}}$ исследованы в~\cite{AHHFG} и~\cite{Ar_smash}. Обе они являются специализациями алгебры ${\mathscr A}_{\omega^\infty}(G)$, определённой выше.

Для произвольной комплексной группы Ли $G$ с алгеброй Ли $\mathfrak{g}$ рассмотрим гомоморфизм
\begin{equation}\label{taudef}
\tau\!:U(\mathfrak{g}) \to \mathscr{A}(G)\!:\langle \sigma(X), f\rangle\!:= [{\widetilde X}f](e) \qquad (X\in U(\mathfrak{g}),\, f \in  \mathcal{O}(G)),
\end{equation}
где ${\widetilde X}$ --- левоинвариантный дифференциальный оператор, соответствующий $X\in U(\mathfrak{g})$.
Следующее предложение в понадобится доказательстве теоремы~\ref{AMheiffsol}, см. \S\,\ref{s:proofUg}.

\begin{pr}\label{evUgAG}
Пусть $\mathcal{C}$ --- класс банаховых алгебр, стабильный относительно перехода к замкнутым подалгебрам и конечным произведениям. Если $G$ односвязна, то гомоморфизм $\widehat U(\mathfrak{g})^\mathcal{C}\to \widehat{\mathscr A}(G)^\mathcal{C}$, полученный применением функтора оболочки к гомоморфизму $\tau$, определённому в~\eqref{taudef}, являeтся изоморфизмом.
\end{pr}
Сначала докажем лемму.
\begin{lm}\label{BwhA}
Пусть $\theta\!:A\to B$ --- эпиморфизм $\ptn$-алгебр.
Если существует непрерывный гомоморфизм  $j\!:B\to \widehat{A}^\mathcal{C}$, такой что $A\to \widehat A^\mathcal{C}$ совпадает с $j\theta$, то $\widehat{\theta}^\mathcal{C}\!:\widehat{A}^\mathcal{C}\to \widehat{B}^\mathcal{C}$ является изоморфизмом $\ptn$-алгебр.
\end{lm}
\begin{proof}
В случае, когда $\mathcal{C}$ --- класс всех банаховых алгебр, утверждение доказано в \cite[Lemma 2.3]{ArAMN}. В общем случае рассуждение такое же, с учётом того, что, как легко видеть, \eqref{AMen} выполнено для всех алгебр, локально содержащихся в~$\mathcal{C}$.
\end{proof}

\begin{proof}[Доказательство предложения~\ref{evUgAG}]
В случае, когда $\mathcal{C}$ --- класс всех банаховых алгебр, в \cite[Proposition~2.1]{ArAMN}  получен
изоморфизм $\widehat U(\mathfrak{g})\cong \widehat{\mathscr A}(G)$. В общем случае рассуждаем аналогичным образом. А именно,
в \cite[Proposition~9.1]{Pir_stbflat} показано, что любой непрерывный гомоморфизм из $U(\mathfrak{g})$ в алгебру Аренса--Майкла пропускается через $\tau$. Так как образ $\tau$ плотен (см. формулу~(42) в \cite{Pir_stbflat} и обсуждение после неё),   из леммы~\ref{BwhA} следует, что $\widehat U(\mathfrak{g})^\mathcal{C}\cong \widehat{\mathscr A}(G)^\mathcal{C}$.
\end{proof}

\begin{rem}
В частных случаях, когда  $\mathcal{C}$ --- класс всех банаховых алгебр и класс банаховых PI-алгебр из предложения~\ref{evUgAG} для односвязной $G$ получаем изоморфизмы $\widehat U(\mathfrak{g})\cong \widehat{\mathscr A}(G)$ и $\widehat U(\mathfrak{g})^{\mathrm{PI}}\cong \widehat{\mathscr A}(G)^{\mathrm{PI}}$. Из второго изоморфизма можно получить следующий факт, упомянутый в  \cite[замечание 6.12]{Ar_smash}: для нильпотентной $\mathfrak{g}$ алгебра $\widehat U(\mathfrak{g})^{\mathrm{PI}}$ топологически изоморфна алгебре формально-радикальных функций  в смысле Доси (последняя определена и исследована в его статьях \cite{Do09C,Do10C,Do10A,Do10B}).

Заметим также, что аналитическая структура $\widehat U(\mathfrak{g})$ для нильпотентной $\mathfrak{g}$ описана в \cite[Theorem 4.3]{AHHFG}, а алгебраическая уточнена в  \cite[Theorem 4.3]{Ar_smash}.
\end{rem}

\subsection*{Гомологические эпиморфизмы и др.}

Напомним  основные определения из гомологической теории локально выпуклых алгебр, подробности  см. в~\cite{X1} и \cite{Pir_qfree}. Мы рассматриваем $\ptn$-алгебры и $\ptn$-модули над ними, т.е. полные локально выпуклые алгебры и модули с совместно непрерывным умножением. Все алгебры предполагаются унитальными.

$A$-$\ptn$-модуль $P$ (левый, правый или бимодуль) называется \emph{проективным}, если для любого допустимого эпиморфизма $A$-$\ptn$-модулей (т.е. обладающего правым обратным непрерывным линейным отображением) существует правый обратный морфизм  $A$-$\ptn$-модулей. Частным случаем  проективного левого модуля является \emph{свободный}, т.е. изоморфный модулю вида $A\ptn E$ для некоторого полного локально выпуклого пространства~$E$.

Если  $M$ и $N$ --- правый  и соответственно левый $A$-$\ptn$-модули, то их \emph{$A$-модульное тензорное произведение} $M\ptens{A}N$ определяется как пополнение факторпространства $M\ptn N$ по замыканию линейной оболочки всех элементов вида 
$$
x\cdot a\otimes y-x\otimes a\cdot y\qquad (x\in M,\,y\in N,\,a\in A).
$$
Цепной комплекс $$\cdots\leftarrow M_n \leftarrow M_{n+1} \leftarrow \cdots$$ $A$-$\ptn$-модулей называется \emph{допустимым}, если он стягиваем в категории топологических векторных пространств, т.е. если существует стягивающая гомотопия, состоящая из непрерывных линейных отображений.

Пусть $R$ --- $\ptn$-алгебра. Напомним, что \emph{$R$-$\ptn$-алгеброй} называется пара $(A,\eta_A)$, где $A$ --- $\ptn$-алгебра, а $\eta_A \!: R \to A$ --- гомоморфизм $\ptn$-алгебр. Заметим, что каждый $A$-$\ptn$-модуль автоматически является $R$-$\ptn$-модулем через функтор ограничения скаляров по $\eta_A$. Обозначим через $(A,R)\lmod$ категорию левых $A$-$\ptn$-модулей, снабжённую  структурой точной категории относительно комплексов, расщепимых непрерывными морфизмами $R$-$\ptn$-модулей; ср. \cite[Appendix, Example 10.1 и 10.3]{Pir_qfree}. В частности, когда $R = \mathbb{C}$, мы получаем стандартное определение допустимой (или $\mathbb{C}$-расщепляемой) последовательности $A$-$\ptn$-модулей, используемое в \cite{X1}. Рассматривая $\ptn$-бимодули над $R$-$\ptn$-алгеброй $A$ (слева) и $S$-$\ptn$-алгеброй $B$ (справа), мы обозначаем соответствующую категорию через $(A,R)\bimod(B,S)$. В случае, когда $R=S=\mathbb{C}$, мы пишем просто $A\bimod B$. Для категорий левых и правых модулей используются обозначения
$A\lmod$ и $\rmod B$ соответственно.

\emph{Проективной (свободной) резольвентой} $M\in A\lmod$ называется допустимый комплекс $0\leftarrow M \xleftarrow{\varepsilon} P_\bullet$, такой, что все модули $P_n$ ($n\ge 0$) проективны  (свободны). Очевидно, всякая свободная резольвента является проективной. Когда мы говорим о проективной резольвенте в $(A,R)\bimod(B,S)$, мы имеем в виду комплекс, состоящий из объектов, проективных в соответствующей точной категории и расщепляющихся в $R\bimod S$.

Пусть $\varphi\!:A\to B$ --- гомоморфизм $\ptn$-алгебр. Он называется \emph{гомологическим эпиморфизмом} (вообще говоря ``сильным'', если нужно отличать от ``слабых'', см. \cite[Definition 3.15]{AP}), если для некоторой (или, эквивалентно, для каждой) проективной резольвенты  $0 \leftarrow B \leftarrow P_\bullet$ в категории левых $A$-$\ptn$-модулей
последовательность
$$
0 \longleftarrow B \longleftarrow B\ptens{A} P_\bullet
$$
допустима. (Это понятие может выступать под разными именами, см. \cite[Remark 3.16]{AP}.)

В доказательстве теоремы~\ref{Assmpr0} нам понадобятся также более общие понятия. 
Напомним, что $R$-$S$-гомоморфизм из $R$-$\mathop{\widehat\otimes}$-алгебры $A$ в $S$-$\mathop{\widehat\otimes}$-алгебру $B$ определяется как пара $(f,g)$, где $f\!: A \to B$ и $g\!: R \to S$ являются гомоморфизмами $\mathop{\widehat\otimes}$-алгебр, такими что $f\eta_A=\eta_Bg$. Далее, $R$-$S$-гомоморфизм $f\!: A \to B$ из $R$-алгебры $A$ в $S$-алгебру $B$ называется \emph{двусторонним относительным гомологическим эпиморфизмом}, если $f$ является эпиморфизмом $\mathop{\widehat\otimes}$-алгебр и при этом модуль $A$ является ациклическим относительно функтора
$$B\ptens{A}(-)\ptens{A}B\!: (A,R){\mbox{-}\!\mathop{\mathsf{mod}}\!\mbox{-}}(A,R)\to (B,S){\mbox{-}\!\mathop{\mathsf{mod}}\!\mbox{-}}(B,S),$$
т.е., этот функтор отображает некоторую (эквивалентно, каждую) проективную резольвенту модуля $A$ в категории $(A,R){\mbox{-}\!\mathop{\mathsf{mod}}\!\mbox{-}}(A,R)$ в последовательность, расщепимую в $S{\mbox{-}\!\mathop{\mathsf{mod}}\!\mbox{-}} S$ \cite[определение 6.2]{Pir_qfree}. 
Мы опускаем определения левых и правых относительных гомологических эпиморфизмов; подробности см. в \cite[определение 6.1]{Pir_qfree} с исправлениями в \cite{Pir_co1,Pir_co2}. Единственное, что мы используем, --- это факт, что двусторонний относительный гомологический эпиморфизм также является левым и правым \cite[предложение 6.5]{Pir_qfree}.

Напомним также, что для $M\in \rmod A$, $N\in A\lmod$ и $n\in\mathbb{Z}_+$ производный функтор ${\mathop{\mathrm{Tor}}\nolimits}_n^{A}(M,N)$ (в классическом смысле) может быть определен как гомология комплекса $M\ptens{A} P_\bullet$, где $P_\bullet$ --- произвольная проективная резольвента $N$, см. \cite[глава III, \S\,4]{X1} или \cite{He81}.

\subsection*{Алгебры степенных рядов}

Для $s\ge 0$ положим
\begin{equation}
 \label{faAsdef}
\mathfrak{A}_s\!:=\Bigl\{a=\sum_{n=0}^\infty  a_n x^n\! :
\|a\|_{r,s}\!:=\sum_{n=0}^\infty |a_n|\frac{r^n}{n!^s}<\infty
\;\forall r>0\Bigr\},
\end{equation}
где $x$ --- формальная переменная.
Для нас важно, что $\mathfrak{A}_s$ является алгеброй Фреше-Аренса--Майкла  \cite[предложение~4]{ArRC} и, более того, $\ptn$-алгеброй Хопфа \cite[Example 2.4]{AHHFG} (операции продолжаются по непрерывности с алгебры многочленов $\mathbb{C}[x]$). Обозначим также через $\mathfrak{A}_\infty$   алгебру $\mathbb{C}[[z]]$ всех формальных степенных рядов.

Связь с алгебрами вида ${\mathscr A}_{\omega^\infty}(\mathbb{C})$  отражена в следующей лемме.

\begin{lm}\label{1dimdsp}
\cite[лемма~6.2]{Ar_smash}
\emph{(A)} Если $\omega(z)=1+|z|$, то ${\mathscr A}_{\omega^\infty}(\mathbb{C})$ топологически изоморфна $\ptn$-алгебре Хопфа $\mathbb{C}[[x]]$ всех формальных рядов от~$x$.

\emph{(B)} Если $s\in [1,\infty)$ и $\omega(z)=\exp(|z|^{1/s})$, то
${\mathscr A}_{\omega^\infty}(\mathbb{C})$ топологически изоморфна $\ptn$-алгебре Хопфа $\mathfrak{A}_{s-1}$.
\end{lm}

Кроме того, следующее предложение существенно используется в дальнейшем.
\begin{pr}\label{AsHomep}
\cite[Proposition 4.4]{Ar_smhe}
Пусть $s\in[0, +\infty]$. Тогда вложение $\mathbb{C}[x]\to \mathfrak{A}_s$ является гомологическим эпиморфимом.
\end{pr}

\subsection*{Смэш-произведения}
Здесь мы коротко напомним необходимые сведения об аналитических смэш-произведениях, впервые рассмотренных в \cite{Pi4}.
Пусть $H$ --- $\ptn$-алгебра Хопфа и пусть $\ptn$-алгебра~$A$  снабжена структурой левого $H$-$\ptn$-модуля. Тогда $A\ptn A$ и $\mathbb{C}$ являются левыми $H$-$\ptn$-модулями относительно внешних умножений, индуцированных коумножением и коединицей. Тогда~$A$ называется (левой)  \emph{$H$-$\ptn$-модульной алгеброй}, если она  снабжена структурой левого $H$-$\ptn$-модуля таким образом, что морфизм умножения и морфизм единицы $\mathbb{C}\to A$ являются морфизмами левых $H$-модулей.  Для всякого  гомоморфизма $\ptn$-алгебр $H\to B$ задано \emph{присоединённое действие} $H$ на $B$ (слева), индуцированное коумножением  и антиподом.

\begin{df}\label{PsmaDef}
Пусть $H$ --- $\ptn$-алгебра Хопфа и $A$ --- $H$-$\ptn$-модульная алгебра.  \emph{Аналитическим  смэш-произведением}  $A\Psmash H$ называется  $\ptn$-алгебра, снабжённая гомоморфизмами $\ptn$-алгебр
$i\!:A\to A\Psmash H$ и $j\!:  H\to A\Psmash H$, такими что выполнены следующие свойства.

(A)~$i$ является гомоморфизмом $H$-$\ptn$-модульных алгебр (относительно
присоединёного действия, порождённого $j$).

(B)~Для любых $\ptn$-алгебры $B$, гомоморфизма $\ptn$-алгебр $\psi\!: H\to B$ и гомоморфизма $H$-$\ptn$-модульных алгебр $\varphi\!: A\to B$ (относительно присоединёного действия, порождённого $\psi$) найдётся единственный гомоморфизм $\ptn$-алгебр $\theta$, такой что диаграмма
\begin{equation*}
   \xymatrix{
 & A\Psmash H \ar@{-->}[dd]^{\theta}& \\
A \ar[ur]^i \ar[dr]_\varphi&&H\ar[ul]_j\ar[dl]^\psi\\
& B &}
\end{equation*}
коммутативна. 
\end{df}

Подробности см. в \cite[\S\,1]{Ar_smash}.

\section{Гомологические эпиморфизмы: формулировка основных результатов}
\label{s:mainfirst}

В этом параграфе содержатся формулировки теорем о гомологических эпиморфизмах. Их доказательства даны в \S\S\,\ref{s:proofmaingen}--\ref{s:proofUg}. 

\subsection*{Алгебры аналитических функционалов}

Сначала сформулируем утверждения об алгебре аналитических функционалов и покажем, что они следуют из более общей теоремы, которую, как мы увидим, можно свести к случаю линейных групп, а затем к случаю односвязных разрешимых.

Основным результатом статьи является следующее утверждение.
\begin{thm}\label{C4main}
Пусть $G$ --- связная комплексная группа Ли. Тогда 
$$
{\mathscr A}(G)\to\widehat{\mathscr A}(G)\quad\text{и}\quad{\mathscr A}(G)\to
\widehat{\mathscr A}(G)^{\mathrm{PI}}
$$
являются гомологическими эпиморфизмами.
\end{thm}

Эта теорема является частным случаем более общей.
Напомним, что  \emph{линеaризатором комплексной группы Ли} называется пересечение ядер всех её конечномерных голоморфных представлений. Все остальные понятия, фигурирующие в следующей теореме о максимальном весе, описаны в \S\,\ref{sec:prel}.

\begin{thm}\label{C4mainexdi}
Пусть $G$ --- связная комплексная группа Ли,  $\Lambda$ --- её линеаризатор, а $E$ и $N$ --- экспоненциальный и нильпотентный радикалы~$G/\Lambda$. Предположим, что  $N'$ ---  нормальная интегральная подгруппа $G/\Lambda$ такая, что $E\subset N'\subset N$, а $\omega_{max}$ --- субмультипликативный вес на~$G/\Lambda$, максимальный среди весов с экспоненциальным искривлением на $N'$.
Тогда ${\mathscr A}(G)\to {\mathscr A}_{\omega_{max}^\infty}(G/\Lambda)$ является гомологическим эпиморфизмом.
\end{thm}

Предполагая теорему~\ref{C4mainexdi} доказанной, мы получаем  теорему~\ref{C4main} следующим образом.
\begin{proof}[Доказательство теоремы~\ref{C4main}]
С одной стороны, так как $G$ связна, то индуцированные факторотображением  $G\to G/\Lambda$ гомоморфизмы
$$
{\widehat{\mathscr A}}(G)\to {\widehat{\mathscr A}}(G/\Lambda)\quad\text{и}\quad
{\widehat{\mathscr A}}(G)^{\mathrm{PI}}\to {\widehat{\mathscr A}}(G/\Lambda)^{\mathrm{PI}}
$$ являются топологическими изоморфизмами
\cite[предложение~5.8]{Ar_smash}. С другой стороны, группа $G/\Lambda$ не только связна, но и линейна, поэтому из \cite[предложение~5.12]{Ar_smash} получаем в случае, если $N'=E$, что ${\mathscr A}_{\omega_{max}^\infty}(G/\Lambda)\cong \widehat{\mathscr A}(G/\Lambda)$, а в случае, если $N'=N$, что ${\mathscr A}_{\omega_{max}^\infty}(G/\Lambda)\cong \widehat{\mathscr A}(G/\Lambda)^{\mathrm{PI}}$.
Таким образом, из теоремы~\ref{C4mainexdi} следует, что ${\mathscr A}(G)\to\widehat{\mathscr A}(G)$ и ${\mathscr A}(G)\to
\widehat{\mathscr A}(G)^{\mathrm{PI}}$ являются гомологическими эпиморфизмами.
\end{proof}

Подробное доказательство теоремы~\ref{C4mainexdi}, которое будет нашей основной целью в дальнейшем, см. в \S\,\ref{s:proofmainexdi}. Здесь же отметим, что оно состоит из трёх шагов:
\begin{itemize}
  \item редукция к случаю группы Штейна,
  \item редукция к случаю линейной группы,
  \item доказательство в случае линейной группы.
\end{itemize}

В случае линейной группы теорема~\ref{C4mainexdi} является следствием другого общего результата --- теоремы о  весе, допускающем разложение, согласованное с итерированным полупрямым произведением.
\begin{df}\label{deflocwe}
Пусть~$\omega$ --- субмультипликативный вес на комплексной группе Ли~$G$, такой что $\omega\ge 1$.  Будем называть~$\omega$ \emph{локализующим}, если эпиморфизм ${\mathscr A}(G)\to{\mathscr A}_{\omega^\infty}(G)$ является гомологическим.
\end{df}

О весах, не являющихся локализующими, см. ниже пример~\ref{nonlovw}.

\begin{thm}\label{C4maingen}
Пусть $G$ ---  комплексная группа Ли и  $\omega$ --- асимптотически симметричный субмультипликативный вес. Предположим, что
$G$ может быть разложена в итерированное полупрямое произведение:
\begin{equation}\label{itdirp0}
 G=(\cdots (F_1 \rtimes F_2)\rtimes\cdots)
\rtimes F_n,
\end{equation}
где $F_1\cong\cdots\cong F_{n-1}\cong\mathbb{C}$, а
$F_n$ линейно комплексно редуктивна. Пусть $G_1\!:=F_1$ и $G_i\!:=G_{i-1}\rtimes F_i$ при $i\le 2$, а $\widetilde\omega_i$ и $\omega_i$ --- ограничения~$\omega$ на~$G_i$ и~$F_i$ соответственно. Предположим также, что имеет место итерированное разложение $\omega$ в следующем смысле:   $\widetilde\omega_1 \simeq  \omega_1$ и для всех $i=1,\ldots,n-1$
\begin{equation}\label{weidec}
\widetilde\omega_{i+1}(gf)\simeq  \widetilde\omega_i(g)\omega_{i+1}(f)\quad\text{на
$G_i\times F_{i+1}$} \quad(\text{$g\in G_i$, $f\in F_{i+1}$}),
\end{equation}
и $\omega_n$ --- максимальный субмультипликативный вес на $F_n$.
Тогда если все веса $\omega_1,\ldots,\omega_n$ --- локализующие, то $\omega$ также является локализующим.
\end{thm}

Доказательство теоремы~\ref{C4maingen} состоит из двух шагов:
\begin{itemize}
  \item редукция к случаю односвязной разрешимой группы (т.е. когда $F_n$ тривиальна),
  \item доказательство в случае односвязной разрешимой группы.
\end{itemize}

Подробности см. в \S\,\ref{s:proofmaingen}.

\subsection*{Универсальные обертывающие алгебры}

Кроме того, в этой статье попутно доказан аналог следующего утверждения об универсальных обертывающих алгебрах.

\begin{thm}\label{AMheiffsol}
\cite{Pi4,Ar_smhe}
Пусть $\mathfrak{g}$ --- конечномерная комплексная алгебра Ли. Оболочка Аренса--Майкла $U(\mathfrak{g})\to \widehat U(\mathfrak{g})$  является гомологическим эпиморфизмом тогда и только тогда, когда~$\mathfrak{g}$ разрешима.
\end{thm}
Необходимость доказана  Пирковским в~\cite{Pi4}, а достаточность автором в~\cite{Ar_smhe}. Отметим, что  в~\cite{ArAMN} для случая, когда $\mathfrak{g}$ нильпотентна,  было дано  доказательство,  основанное на методах, отличных от использованных в~\cite{Ar_smhe}. Здесь мы докажем аналогичное утверждение для оболочки в классе  банаховых PI-алгебр.

\begin{thm}\label{AMheifnew}
Пусть $\mathfrak{g}$ --- конечномерная комплексная алгебра Ли. Оболочка в классе  банаховых PI-алгебр
$U(\mathfrak{g})\to \widehat U(\mathfrak{g})^{\mathrm{PI}}$   является гомологическим эпиморфизмом
тогда и только тогда, когда~$\mathfrak{g}$ разрешима.
\end{thm}

Доказательство этой теоремы см. в \S\,\ref{s:proofUg}.

\section{Доказательство теоремы о  весе, допускающем разложение}
\label{s:proofmaingen}

В формулировке теоремы~\ref{C4maingen} (о  весе, допускающем разложение, согласованное с итерированным полупрямым произведением) последняя подгруппа $F_n$, участвующая в произведении, предполагается линейно комплексно редуктивной. Мы увидим, что соответствующая алгебра ${\mathscr A}(F_n)$ является гомологически тривиальной. Первый шаг доказательства теоремы --- редукция к случаю односвязной разрешимой группы --- опирается на результаты о гомологически тривиальных алгебрах. Мы обсудим их сначала в абсолютном, а затем в относительном вариантах.

\subsection*{Случай линейно комплексно редуктивной группы}

Напомним, что комплексная группа Ли~$G$  называется \emph{линейно комплексно редуктивной}, если существует компактная действительная группа Ли, такая что $G$ является её универсальной комплексификацией  \cite[Definition~15.2.7]{HiNe}. (Связность здесь не требуется, но очевидно, что $G$ имеет лишь конечное число связных компонент.)

Напомним также, что унитальная $\ptn$-алгебра $A$  называется  \emph{гомологически тривиальной} (или \emph{стягиваемой}), если $A$ является проективным $A$-$\ptn$-бимодулем \cite{X1}.

Отметим, что следующая теорема  является аналогом такого результата Тейлора: если $K$ --- компактная действительная группа Ли, то алгебра распределений ${\mathscr E}'(K)$ гомологически тривиальна  \cite[Proposition~7.3]{T2}. Указанная теорема Тейлора также является важной частью нашего доказательства.

\begin{thm}\label{Aredhtr}
Пусть $G$ --- линейно комплексно редуктивная группа Ли. Тогда ${\mathscr A}(G)$ гомологически тривиальна.
\end{thm}

Нам понадобится несколько лемм.

\begin{lm}\label{derhtr}
Пусть $\sigma\!:A\to B$ --- гомоморфизм $\ptn$-алгебр с плотным образом. Если~$A$ гомологически тривиальна, то такова же и~$B$.
\end{lm}
\begin{proof}
В силу \cite[предложение~IV.5.7]{X1}  $\ptn$-алгебра $A$ гомологически тривиальна тогда и только тогда $A\ptn A$ содержит диагональ, т.е. элемент~$d$, для которого выполнены два условия:

(1)~$a\cdot d = d\cdot a$ каждого $a\in A$;

(2)~умножение $A\ptn A \to A$ отображает~$d$ в единицу~$A$.

Если $d$ --- диагональ в $A\ptn A$, то из того,что образ $\sigma$ плотен, получаем, что $(\sigma\otimes\sigma)(d)$ является диагональю в $B\ptn B$.
\end{proof}

Ниже мы обозначаем через $C^\infty(G)$ локально выпуклое пространство комплекснозначных бесконечнодифференцируемых функций на действительной группе Ли~$G$, а через ${\mathscr E}'(G)$ ---  пространство распределений с компактным носителем на $G$, т.е. сильное двойственное пространство к $C^\infty(G)$. Хорошо известно, что ${\mathscr E}'(G)$ является $\ptn$-алгеброй относительно свёртки. (Оно даже является $\ptn$-алгеброй Хопфа, но нам это здесь не нужно.) Мы отождествляем ${\mathscr E}'(G)\ptn {\mathscr E}'(G)$ с сильным двойственным пространством к $C^\infty(G\times G)$.

\begin{lm}\label{diagcogr}
Пусть $K$ --- компактная действительная группа Ли.
Тогда  функционал
$$
d\!:C^\infty(K\times K)\to \mathbb{C}\!:f\mapsto  \int_K f(k,k^{-1})\,dk
$$
где интегрирование производится по нормированной мере Хаара, является диагональю в ${\mathscr E}'(K)\ptn {\mathscr E}'(K)$.
\end{lm}
Доказательство принадлежит Дж.~Тейлору (см. \cite[Proposition~7.3]{T2}) и мы опускаем его.

Следующая лемма доказана Литвиновым.
\begin{lm}\label{litvden}
\cite[предложение~4]{Li70}
Пусть $G$  --- комплексная группа Ли, обладающая действительной формой $G_\mathbb{R}$. Тогда гомоморфизм групп $G_\mathbb{R}\to G$  порождает гомоморфизм $\ptn$-алгебр 
$$
\sigma\!:{\mathscr E}'(G_\mathbb{R})\to {\mathscr A}(G)\!:\langle \sigma(\alpha), f \rangle = \langle \alpha, f|_{G_\mathbb{R}}\rangle \qquad (\alpha \in {\mathscr E}'(G_\mathbb{R}),\,
f\in \mathcal{O}(G)).
$$ 
Более того, этот гомоморфизм имеет плотный образ.
\end{lm}
Заметим, в частности, что плотность образа следует непосредственно из условий Коши--Римана.

\begin{proof}[Доказательство теоремы~\ref{Aredhtr}]
Пусть $K$ --- компактная действительная группа Ли, такая что $G$ является её универсальной комплексификацией. В силу компактности~$K$ алгебра ${\mathscr E}'(K)$  гомологически тривиальна  \cite[Proposition~7.3]{T2}, ср. с леммой~\ref{diagcogr}. Поскольку $K$ является действительной формой $G$, утверждение следует из лемм~\ref{derhtr} и~\ref{litvden}.
\end{proof}

\begin{rem}\label{expdia}
В случае, когда линейно комплексно редуктивная группа Ли $G$ связна,
теорема~\ref{Aredhtr} может быть доказана с использованием другой явной конструкции диагонали в ${\mathscr A}(G)\ptn{\mathscr A}(G)$, отличной от описанной в лемме~\ref{diagcogr}. (Не исключено, что это возможно и в общем случае.)

В самом деле, пусть $K$ --- компактная действительная группа Ли, такая что $G$ является её универсальной комплексификацией. Обозначим через~$\Sigma$ множество доминантных аналитически интегральных весов (мы следуем терминологии из~\cite{Sa02}). Тогда в силу теоремы о старших весах множество неприводимых унитарных представлений~$K$ находится во взаимно-однозначном соответствии с~$\Sigma$. Тем самым преобразование Фурье является топологическим изоморфизмом между $\mathcal{O}(G)$ и пространством
$$
{\mathscr A}_{exp}(\Sigma)\!:= \Bigl\{ s=(s_\sigma;\,\sigma\in \Sigma)\!:\;\forall
r>0\; \|s\|_r\!:=e^{r\|\sigma\|} \|s_\sigma\|_{HS}<\infty \Bigr\},
$$
где $s_\sigma$ --- оператор в пространстве соответствующего представления, $\|\cdot\|_{HS}$ --- норма Гильберта-Шмидта, а $\|\cdot\|$ --- какая-нибудь норма на пространстве линейных функционалов на алгебре Ли, см. \cite[формулы (19)--(22)]{Sa02}.  Следовательно, мы получаем  изоморфизм $\ptn$-алгебр
между ${\mathscr A}(G)$ и
$$
\mathcal{O}_{exp}(\Sigma)\!:= \Bigl\{ t=(t_\sigma;\sigma\in \Sigma)\!:\;\exists
C,\,r>0\; \|t_\sigma\|_{HS}<C\,e^{r\|\sigma\|}  \Bigr\}
$$
(здесь $\mathcal{O}_{exp}(\Sigma)$ есть $\ptn$-алгебра относительно операторного умножения в каждом слагаемом и топологии индуктивного предела). Рассмотрим последовательность (индексированных $\Sigma\times\Sigma$) операторов, принимающую
значение
$$
\sum_{ij}\frac1{\dim \sigma}\,e^\sigma_{ij}\otimes e^\sigma_{ji}
$$
в $(\sigma,\sigma)$ и $0$ для остальных индексов (здесь $e^\sigma_{ij}$ --- матричная единица в пространстве представления $\sigma$). Легко видеть, что мы получили диагональ в $\mathcal{O}_{exp}(\Sigma)\ptn\mathcal{O}_{exp}(\Sigma)$.
\end{rem}

\begin{lm}\label{hoholtr}
Пусть $\varphi\!:A\to B$ --- эпиморфизм $\ptn$-алгебр. Если  $A$ гомологически тривиальна, то $\varphi$ является гомологическим.
\end{lm}
\begin{proof}
Так как $A$ гомологически тривиальна, все модули из $A\lmod$ проективны, в том числе  $B$. Тем самым $0 \longleftarrow B \longleftarrow B \longleftarrow 0$ является проективной резольвентой в $A\lmod$. Применяя функтор $B\ptens{A}(-)$ получаем последовательность
$$0 \longleftarrow B\ptens{A} B \longleftarrow B\ptens{A} B \longleftarrow 0,$$ которая очевидно допустима. Таким образом, $\varphi$ --- гомологический.
\end{proof}

Следующий результат понадобится в доказательстве теоремы~\ref{C4maingen}. Кроме того, если~$\omega$ удовлетворяет условию максимальности, указанному в формулировке теоремы~\ref{C4mainexdi},
с $N' = E$ или $N' = N$, мы получаем утверждение теоремы~\ref{C4main} для линейно комплексно редуктивной группы.

\begin{pr}\label{coredhe}
Пусть~$\omega$ --- субмультипликативный вес на линейно комплексно редуктивной группе Ли~$G$, такой что $\omega\ge 1$. Тогда $\omega$ является локализующим.
\end{pr}
\begin{proof}
В силу теоремы~\ref{Aredhtr} $\ptn$-алгебра ${\mathscr A}(G)$ гомологически тривиальна. Остаётся воспользоваться леммой~\ref{hoholtr}.
\end{proof}

\subsection*{Относительно гомологически тривиальные алгебры}

Теоремы~\ref{Aredhtr} недостаточно для достижения нашей ближайшей цели ---  редукции к случаю односвязной разрешимой группы. Чтобы получить её обобщение --- теорему~\ref{redrelhotr}, мы используем относительную гомологическую теорию.

\begin{df}
Будем говорить, что унитальная $R$-$\ptn$-алгебра $A$ \emph{относительно гомологически тривиальна},  если $A$ является проективным модулем в относительной категории $(A,R)\bimod(A,R)$.
\end{df}

\begin{lm}\label{reldia}
$R$-$\ptn$-алгебра $A$ относительно  гомологически тривиальна тогда и только тогда, когда $A\ptens{R} A$ содержит \emph{относительную диагональ}, т.е. элемент~$d$, для которого выполнены два условия:

\emph{(1)}~$a\cdot d = d\cdot a$ каждого $a\in A$;

\emph{(2)}~умножение $A\ptens{R}A \to A$ отображает~$d$ в единицу~$A$.
\end{lm}
Доказательство в точности то же, что и для случая $R=\mathbb{C}$, см. \cite[предложение~IV.5.7]{X1}.

\begin{thm}\label{redrelhotr}
Пусть комплексная группа Ли имеет вид $G=N\rtimes L$, где $L$ линейно комплексно редуктивна.  Тогда ${\mathscr A}(N)$-$\ptn$-алгебра  ${\mathscr A}(G)$ относительно гомологически тривиальна.
\end{thm}

Для доказательства потребуются вспомогательные утверждения. 

\begin{lm}\label{deltden}
Линейная оболочка множества дельта-функций  на комплексном многообразии плотна в пространстве аналитических функционалов.
\end{lm}
\begin{proof}
B силу теоремы Хана-Банаха достаточно убедиться, что всякий непрерывный линейный функционал на пространстве аналитических функционалов, значения которого равны $0$ на всех дельта-функциях, тождественно равен $0$. Последнее выполнено, поскольку непрерывный линейный функциoнал определяется голоморфной функцией на  многообразии ввиду рефлексивности пространства голоморфных функций.
\end{proof}

Рассмотрим более общую ситуацию. Пусть $N$ --- замкнутая комплексная подгруппа в комплексной группе Ли~$G$. Обозначим  подпространство в $\mathcal{O}(G\times G)$, состоящее из функций, удовлетворяющих условию 
\begin{equation}\label{balfun}
f(g_1n,g_2)=f(g_1,ng_2)\qquad (n\in N,\,g_1,g_2\in G),
\end{equation}
через $\mathcal{O}_{bal}(G\times G)$.
Легко видеть, что оно замкнуто в $\mathcal{O}(G\times G)$.

Заметим, что $\mathcal{O}(G\times G)$ является ${\mathscr A}(G)$-$\ptn$-бимодулем относительно действий
$$
\alpha\bullet f = f \ast \widetilde\alpha \quad\text{и} \quad f \bullet  \alpha = \widetilde \alpha \ast f \qquad (\alpha\in {\mathscr A}(G),\,f\in \mathcal{O}(G\times G));
$$ 
здесь тильда обозначает стандартный антипод в ${\mathscr A}(G)$. Кроме того, мы рассматриваем ${\mathscr A}(G\times G)$ как ${\mathscr A}(G)$-$\ptn$-бимодуль, отождествляя его с   ${\mathscr A}(G)\ptn {\mathscr A}(G)$.

Исходя из стандартных формул для свертки функции и дельта-функции (см., например, \cite[\S\,3]{Ak08}), получаем, что
\begin{equation}\label{degfunc}
(f \bullet \delta_{g})(h_1,h_2)=f(gh_1,h_2)\quad\text{и}\quad(\delta_{g}\bullet f)(h_1,h_2)=f(h_1,h_2g)
\end{equation}
для всех $g\in G$ и $f\in\mathcal{O}(G\times G)$. Отсюда сразу видно, что $\mathcal{O}_{bal}(G\times G)$ инвариантно относительно действия $\delta_{g}$ слева и справа. Из леммы~\ref{deltden} о плотности следует, что $\mathcal{O}_{bal}(G\times G)$ является подмодулем в ${\mathscr A}(G)$-$\ptn$-бимодулe $\mathcal{O}(G\times G)$. Далее мы предполагаем, что $\mathcal{O}_{bal}(G\times G)'$ снабжено стандартной структурой двойственного локально выпуклого ${\mathscr A}(G)$-бимодуля, т.е.
\begin{equation}\label{degfunc2}
\langle \alpha\bullet \eta ,f \rangle = \langle \eta ,f\bullet \alpha   \rangle\quad\text{и} \quad \langle\eta \bullet \alpha  ,f \rangle = \langle \eta ,\alpha\bullet f \rangle,
\end{equation}
где $\alpha\in {\mathscr A}(G)$, $\eta\in \mathcal{O}_{bal}(G\times G)'$ и $f\in \mathcal{O}_{bal}(G\times G)$.

\begin{pr}\label{norbal2}
Изоморфизм локально выпуклых пространств 
$$
\upsilon\!:{\mathscr A}(G)\ptn{\mathscr A}(G)\to \mathcal{O}(G\times G)'\!:\delta_{h_1}\otimes \delta_{h_2}\mapsto (f\mapsto f(h_1,h_2))
$$
порождает линейное отображение 
$$
\Upsilon\!:{\mathscr A}(G)\ptens{{\mathscr A}(N)} {\mathscr A}(G)\to \mathcal{O}_{bal}(G\times G)',
$$
которое является топологическим изоморфизмом ${\mathscr A}(G)$-бимодулeй.
\end{pr}

\begin{proof}
По определению, 
${\mathscr A}(G)\ptens{{\mathscr A}(N)} {\mathscr A}(G)$ является пополнением факторпространства $({\mathscr A}(G)\ptn {\mathscr A}(G))/F $, где 
$F$ --- замыкание в ${\mathscr A}(G)\ptn {\mathscr A}(G)$ линейной оболочки элементов вида 
$$
(\alpha\ast\eta)\otimes \beta- \alpha\otimes (\eta\ast\beta) \qquad (\alpha, \beta\in {\mathscr A}(G),\,\eta\in {\mathscr A}(N))
$$ 
\cite[теорема II.4.11]{X1}.
Из~\eqref{balfun} следует, что если  $n\in N$, $h_1,h_2\in G$, то для всех $f\in\mathcal{O}_{bal}(G\times G)$ выполнено равенство
\begin{equation}\label{FannOb}
\langle (\delta_{h_1}\ast \delta_n)\otimes \delta_{h_2}, f\rangle=f(h_1n,h_2)=f(h_1,nh_2)=
\langle \delta_{h_1} \otimes  (\delta_n\ast \delta_{h_2}), f\rangle.
\end{equation}
Тем самым $\upsilon((\delta_{h_1}\ast\delta_n)\otimes  \delta_{h_2}- \delta_{h_1}\otimes (\delta_n\ast \delta_{h_2}))$ содержится в $\mathcal{O}_{bal}(G\times G)^\perp$,
где $\perp$ обозначает аннулятор. Тогда из леммы~\ref{deltden} о плотности 
следует, что  
$$
\upsilon((\alpha\ast \eta)\otimes \beta- \alpha\otimes (\eta\ast\beta))\in \mathcal{O}_{bal}(G\times G)^\perp
$$ 
для всех $\alpha,\beta\in{\mathscr A}(G)$ и $\eta\in{\mathscr A}(N)$. Следовательно, мы получаем непрерывное линейное отображение 
$$
{\mathscr A}(G)\ptens{{\mathscr A}(N)} {\mathscr A}(G)\to \mathcal{O}(G\times G)'/\mathcal{O}_{bal}(G\times G)^\perp.
$$

Заметим, что всякое ядерное пространство является пространством Шварца \cite[Corollary 28.5]{MV}. Тем самым $\mathcal{O}(G\times G)$ есть пространство Фреше--Шварца, а значит к нему применимо \cite[Corollary 26.25]{MV}, которое обеспечивает канонический топологический изоморфизм $\mathcal{O}_{bal}(G\times G)'\to \mathcal{O}(G\times G)'/\mathcal{O}_{bal}(G\times G)^\perp$. Итак, линейное отображение $\Upsilon$ корректно определено, непрерывно и открыто.

Чтобы показать, что $\Upsilon$ является топологическим изоморфизмом, достаточно установить равенство $\upsilon(F) = \mathcal{O}_{bal}(G\times G)^\perp$.
Действительно, из \eqref{FannOb} вытекает, что $\upsilon(F)^\perp=\mathcal{O}_{bal}(G\times G)$. Поэтому из теоремы о двойном аннуляторе следует, что
$\mathcal{O}_{bal}(G\times G)'$ есть замыкание $\upsilon(F)$ в слабой топологии $\sigma(\mathcal{O}(G\times G)',\mathcal{O}(G\times G))$
пространства $\mathcal{O}(G\times G)'$. Остаётся воспользоваться рефлексивностью пространства $\mathcal{O}(G\times G)$ и очевидной замкнутостью $\upsilon(F)$ в сильной топологии на $\mathcal{O}(G\times G)'$.

Для завершения доказательства необходимо проверить, что $\Upsilon$  есть морфизм ${\mathscr A}(G)$-$\ptn$-бимодулeй. В силу леммы~\ref{deltden} достаточно проверить, что
$$
\Upsilon(\delta_g\ast\delta_{h_1}\otimes \delta_{h_2})=\delta_g\bullet\Upsilon(\delta_{h_1}\otimes \delta_{h_2})\quad\text{и}\quad
\Upsilon(\delta_{h_1}\otimes \delta_{h_2}\ast\delta_g)=\Upsilon(\delta_{h_1}\otimes \delta_{h_2})\bullet\delta_g
$$
для всех $g,h_1,h_2\in G$.

Пусть $f\in\mathcal{O}_{bal}(G\times G)$. Тогда из~\eqref{degfunc} следует, что
$$
\langle\Upsilon(\delta_g\ast\delta_{h_1}\otimes \delta_{h_2}),f\rangle=f(gh_1,h_2)=(f \bullet \delta_{g})(h_1,h_2)=\langle\delta_g\bullet\Upsilon(\delta_{h_1}\otimes \delta_{h_2}), f\rangle.
$$
Первая формула доказана. Вторая проверяется аналогично.
\end{proof}

Теперь мы можем доказать теорему об относительной гомологической тривиальности.

\begin{proof}[Доказательство теоремы~\ref{redrelhotr}]
Согласно лемме~\ref{reldia} достаточно показать, что модуль ${\mathscr A}(G)\ptens{{\mathscr A}(N)} {\mathscr A}(G)$ содержит относительную диагональ.

В силу предложения~\ref{norbal2} элементы модуля ${\mathscr A}(G)\ptens{{\mathscr A}(N)} {\mathscr A}(G)$ можно отождествить с функционалами на $\mathcal{O}_{bal}(G\times G)$. Пусть $K$ --- компактная подгруппа $L$,  для которой последняя является универсальной комплексификацией. Проверим, что  функционал
\begin{equation}\label{forndia}
d\!:\mathcal{O}_{bal}(G\times G)\to \mathbb{C}\!:f\mapsto  \int_K f(k,k^{-1})\,dk,
\end{equation}
где интегрирование производится по нормированной мере Хаара,
является относительной диагональю. (Мы используем ту же формулу, что и в лемме~\ref{diagcogr}.) Достаточно проверить  условия~(1) и~(2) леммы~\ref{reldia}.

(1)~Напомним, что  согласно лемме~\ref{deltden} линейная оболочка множества дельта-функций плотна в ${\mathscr A}(G)$. Поэтому условие $\mu\bullet d=d\bullet\mu$ достаточно проверить только в случае, когда $\mu=\delta_g$, $g\in G$. Так как~$G$ является полупрямым произведением, то будем проверять отдельно для элементов $N$ и~$L$.

Пусть $n\in N$. Из~\eqref{degfunc2} и \eqref{degfunc} следует, что
$$
 \langle\delta_n\bullet d,f\rangle=\int_K f(nk,k^{-1})\,dk\quad\text{и}\quad \langle d\bullet \delta_n,f\rangle=\int_K f(k,k^{-1}n)\,dk
$$
для любой функции $f\in\mathcal{O}_{bal}(G\times G)$.
Так как для каждого $k\in K$ имеем $nk=kn_1$, где $n_1=k^{-1}nk \in N$, то из условия $f\in \mathcal{O}_{bal}(G\times G)$, следует, что 
$$
f(nk,k^{-1})=f(kn_1,k^{-1})=f(k,n_1k^{-1})=f(k,k^{-1}n).
$$ 
Тем самым $\delta_n\bullet d=d\bullet \delta_n$.

Рассмотрим функционал $\widetilde d$ на $\mathcal{O}(L\times L)$, заданный той же формулой, что и в~\eqref{forndia}. В силу леммы~\ref{diagcogr} та же формула задаёт также диагональ в ${\mathscr E}'(K)\ptn {\mathscr E}'(K)$. Из леммы~\ref{litvden} следует, что вложение $K\to L$  порождает гомоморфизм $\ptn$-алгебр $\sigma\!:{\mathscr E}'(K)\to {\mathscr A}(L)$, который имеет плотный образ. Следовательно,  диагональ в ${\mathscr E}'(K)\ptn {\mathscr E}'(K)$ отображается в диагональ в ${\mathscr A}(L)\ptn{\mathscr A}(L)$, ср.  лемму~\ref{derhtr}. Тем самым $\widetilde d$ есть диагональ в ${\mathscr A}(L)\ptn{\mathscr A}(L)$. В частности, это означает, что $\delta_l\bullet \widetilde d=\widetilde d\bullet \delta_l$ для каждого $l\in L$. Так как естественное отображение ${\mathscr A}(L)\ptn{\mathscr A}(L)\to {\mathscr A}(G)\ptens{{\mathscr A}(N)} {\mathscr A}(G)$ является морфизмом ${\mathscr A}(L)$-бимодулей, получаем, что $\delta_l\bullet d=d\bullet \delta_l$ для каждого $l\in L$.

(2)~Рассмотрим коумножение на $\mathcal{O}(G)$, т.е. отображение $\Delta(h)(g_1,g_2)=h(g_1g_2)$, $g_1,g_2\in G$ (см. начало \S\,\ref{sec:prel}). Тогда умножение $m\!:{\mathscr A}(G)\ptens{{\mathscr A}(N)} {\mathscr A}(G)\to {\mathscr A}(G)$ имеет вид
$$
\langle m(\mu\otimes\eta), h\rangle=\langle \mu\otimes\eta,\Delta(h)\rangle\qquad (\mu,\eta \in {\mathscr A}(G),\, h\in \mathcal{O}(G)).
$$
Следовательно,
$$
\langle m(d),h\rangle=\langle d,\Delta(h)\rangle=\int_K h(kk^{-1})\,dk=h(e)=\langle\delta_e, h\rangle
$$
для всех $h$. Так как $\delta_e$ является единицей в ${\mathscr A}(G)$, условие~(2) леммы~\ref{reldia} выполнено.
\end{proof}

\subsection*{Редукция к случаю односвязной разрешимой группы}

Мы используем гомоморфизмы между аналитическими смэш-произведениями. Существование таковых обеспечивается следующей леммой.

\begin{lm}\label{homsmpr}
\cite[лемма 3.10]{Ar_smash}
Пусть $H$ и $K$ --- $\mathbin{\widehat{\otimes}}$-алгебры Хопфа, $R$ --- $H$-$\mathbin{\widehat{\otimes}}$-модульная алгебра, а $S$ --- $K$-$\mathbin{\widehat{\otimes}}$-модульная алгебра. Если $\beta\!:H\to K$ ---  гомоморфизм $\mathbin{\widehat{\otimes}}$-алгебр Хопфа  и  $\alpha\!: R\to S$ --- гомоморфизм $\mathbin{\widehat{\otimes}}$-алгебр, который также является морфизмом $H$-$\mathbin{\widehat{\otimes}}$-модулей, то формула
\begin{equation}\label{morsmpr}
\alpha\mathop{\widehat{\#}}\beta\!:R\mathop{\widehat{\#}}H\to S\mathop{\widehat{\#}} K\!:r\otimes
h\mapsto \alpha(r)\otimes \beta(h)
\end{equation}
задаёт  гомоморфизм  $\mathbin{\widehat{\otimes}}$-алгебр.
\end{lm}

Теперь сформулируем ключевое утверждение об относительно гомологически тривиальных алгебрах (оно параллельно теореме~\ref{Assmpr} ниже).

\begin{thm}\label{Assmpr0}
Пусть выполнены условия леммы~\ref{homsmpr} и, кроме того,

\emph{(1)}~$\alpha$ и $\beta$ гомологические эпиморфизмы и имеют плотные образы,

\emph{(2)}~$R\Psmash H$ является относительно гомологически тривиальной $R$-$\ptn$-алгеброй.

Тогда $\alpha\mathop{\widehat{\#}}\beta\!:R\Psmash H\to S\Psmash K$ также
является гомологическим эпиморфизмом и имеет плотный образ.
\end{thm}

Для доказательства нам понадобятся два свойства относительно гомологически тривиальных алгебр, которые позволяют использовать результаты Пирковского из \cite{Pir_qfree} и их модификации из \cite{Ar_smhe}.

Напомним, что $R$-$\ptn$-алгебра $A$ называется \emph{относительно квазисвободной}, если всякое её допустимое сингулярное $R$-расширение расщепимо, см. \cite[Definition 7.1]{Pir_qfree} и там же обсуждение этого понятия.

\begin{pr}\label{hotrqufr}
Пусть $R$ --- $\ptn$-алгебра. Тогда всякая относительно гомологически тривиальная $R$-$\ptn$-алгебра относительно квазисвободна.
\end{pr}
\begin{proof}
Пусть $A$ --- $R$-$\ptn$-алгебра и пусть $\mathop{\mathrm{db}}_R$ обозначает относительную биразмерность Хохшильда. Как видно из определений, $A$ является относительно гомологически тривиальной тогда и только тогда, когда $\mathop{\mathrm{db}}_R A = 0 $ и относительно квазисвободной  тогда и только тогда, когда $\mathop{\mathrm{db}}_R A \le 1$ \cite[Proposition 7.1]{Pir_qfree}. Ясно, что первое условие
влечёт второе.

\end{proof}

Следуя \cite[Definition 8.3]{Pir_qfree}, мы говорим, что проективный $(A, R)$-$\ptn$-бимодуль $P$ удовлетворяет \emph{условию конечности $(f_2)$}, если он является ретрактом бимодуля вида  $A\ptens{R}R^n\ptens{R} A$ для некоторого $n\in\mathbb{N}$. (Альтернативное условие $(f_1)$ мы не используем.)
Далее, $R$-$\ptn$-алгебра $A$ называется  $R$-$\ptn$-алгеброй \emph{$(f_2)$-конечного  типа}, если $A$ допускает в $(A, R)\bimod(A, R)$ проективную резольвенту конечной длины, состоящую из  бимодулей, удовлетворяющих условию конечности $(f_2)$ \cite[Definition 8.4]{Pir_qfree}.

\begin{pr}\label{hotrf2}
Всякая относительно гомологически тривиальная $R$-$\ptn$-алгебра является алгеброй $(f_2)$-конечного типа.
\end{pr}
\begin{proof}
Пусть $A$ --- относительно гомологически тривиальная $R$-$\ptn$-алгебра. Тогда  она является ретрактом $(A, R)$-$\ptn$-бимодуля $A\ptens{R}A$, который очевидно изоморфен бимодулю $A\ptens{R}R\ptens{R}A$, а последний удовлетворяет условию конечности $(f_2)$. Тем самым $A$ имеет в $(A, R)\bimod(A, R)$ резольвенту длины $0$ и при этом выполнено  условие конечности $(f_2)$.
\end{proof}

Кроме того, в следующем ниже доказательстве теоремы~\ref{Assmpr0} используются понятия гомоморфизма, удовлетворяющего свойству (UDE), и относительного гомологического эпиморфизма. 
По определению, непрерывный гомоморфизм $\varphi\!:A\to B$ между $\mathop{\widehat\otimes}$-алгебрами удовлетворяет  \emph{свойству (UDE)}, если каждое непрерывное дифференцирование из $A$ в произвольный $B$-$\mathop{\widehat\otimes}$-бимодуль $X$ единственным образом продолжается до непрерывного дифференцирования из $B$ в $X$ \cite[Definition 3.1]{Ar_smhe}. Определение второго понятия см. в \S~\ref{sec:prel}.

\begin{proof}[Доказательство теоремы~\ref{Assmpr0}]
Рассуждения близки к приведённым в доказательстве теоремы~4.5  из \cite{Ar_smhe}, но требуются некоторые модификации.

Положим $A=R\Smash H$, $B=S\Psmash K$, $g=\alpha$ и $f=\alpha\mathop{\widehat{\#}}\beta$. Достаточно проверить выполнение условий теоремы~1.1  из \cite{Ar_smhe}, а именно, что

(1)~$g$ является гомологическим эпиморфизмом;

(2)~$R$-$S$-гомоморфизм $(f, g)$ является левым или правым относительным гомологическим эпиморфизмом;

(3)~$A$ проективен в $R\lmod$, а $B$ проективен в $S\lmod$;

(4)~$A$ является $R$-$\mathbin{\widehat{\otimes}}$-алгеброй $(f_2)$-конечного типа.

Условие~(1) входит в формулировку теоремы. Так как двусторонний относительный гомологический эпиморфизм также является левым и правым \cite[Proposition~6.5]{Pir_qfree}, из   \cite[Theorem~4.6]{Ar_smhe} следует, что для проверки условия~(2) достаточно установить следующие факты: 

(A)~$g\!: R \to S$ имеет плотный образ;

(B)~$f\!: A \to B$ удовлетворяет свойству (UDE) и имеет плотный образ;

(C)~$A$ относительно квазисвободна над $R$.

Действительно, плотность образов $f$ и $g$ следует из условий теоремы.
Далее, поскольку $\alpha$ и $\beta$ являются гомологическими эпиморфизмами, свойство (UDE) выполнено для них по теореме~4.11  из \cite{Ar_smhe}. Более того, \cite[теорема~3.6]{Ar_smhe} влечёт, что $f=\alpha\mathop{\widehat{\#}}\beta$ также удовлетворяет свойству (UDE). 
Из  предложения~\ref{hotrqufr} следует, что $A$ относительно квазисвободна над $R$. Итак, условия (A), (B) и~(C) выполнены, что доказывает условие~(2).

Для проверки условия~(3) заметим, что $S\Psmash K$ изоморфен $S\mathbin{\widehat{\otimes}} K$ как $S$-$\mathbin{\widehat{\otimes}}$-модуль, и, следовательно, он проективен, будучи свободным модулем. Аналогично, $R\Smash H$ является проективным $R$-модулем.

Наконец, условие~(4) выполняется в силу предложения~\ref{hotrf2}.
\end{proof}

Предложение~\ref{coredhe} в совокупности с теоремами~\ref{redrelhotr} и~\ref{Assmpr0} позволяет  произвести редукцию теоремы~\ref{C4maingen} к случаю односвязной разрешимой группы. Подробности см. ниже. Сам односвязный разрешимый случай  обсуждается в следующем параграфе.

\subsection*{Доказательство в случае односвязной разрешимой группы}

Случай  односвязной разрешимой группы фактически исследован в \cite{Ar_smhe}, но для алгебр Ли и в несколько меньшей общности. Поэтому требуются уточнения.

Итак, мы сначала докажем частный случай теоремы~\ref{C4maingen}.

\begin{thm}\label{solvmain}
Пусть выполнены условия теоремы~\ref{C4maingen}. Предположим дополнительно, что $F_n$ тривиальна. Тогда $\omega$ является локализующим.
\end{thm}
Заметим, что в этом случае $G$ разрешима и  односвязна. Для доказательства нам понадобится следующая теорема об универсальных обертывающих, которая является простым обобщением утверждения о достаточности в теореме~\ref{AMheiffsol} (которая есть теорема~4.3 из \cite{Ar_smhe}). Причина заключается в том, что в последней используется только существование разложения~\eqref{expfsmp}.

В формулировке теоремы мы используем следующее понятие. Пусть $H$ и $K$ --- $\ptn$-алгебры Хопфa, каждая из которых представлена как итерированное аналитическое смэш-произведение, т.е. 
$H\cong(\cdots (H_1 \mathop{\widehat{\#}}H_2)\mathop{\widehat{\#}}\cdots)\mathop{\widehat{\#}}H_n$,
где $R_1=H_1$, $R_i=R_{i-1}\Psmash H_{i}$ и $R_{i-1}$ является $H_{i}$-$\ptn$-модульной биалгеброй для каждого $i=1,\ldots n-1$, и аналогично для $K\cong(\cdots (K_1 \mathop{\widehat{\#}}K_2)\mathop{\widehat{\#}}\cdots)\mathop{\widehat{\#}}K_n$,
где $S_1=K_1$, $S_i=S_{i-1}\Psmash K_{i}$ и $S_{i-1}$ является $K_{i}$-$\ptn$-модульной биалгеброй для каждого $i=1,\ldots n-1$.  Будем говорить, что гомоморфизм $\ptn$-алгебр $H\to K$ \emph{согласован} с этими итерированными разложениями, если на каждом шаге итерации мы имеем гомоморфизм  $\mathbin{\widehat{\otimes}}$-алгебр Хопфа $H_k\to K_k$ и гомоморфизм $\mathbin{\widehat{\otimes}}$-алгебр $R_{k-1}\to S_{k-1}$, который является морфизмом $H_k$-модулей, ср. лемму~\ref{homsmpr}.

\begin{thm}\label{AMede2}
Пусть $\mathfrak{g}$ — конечномерная разрешимая комплексная алгебра Ли и
\begin{equation}\label{fgdec}
\mathfrak{g}=((\cdots (\mathfrak{f}_1 \rtimes \mathfrak{f}_2)\rtimes\cdots)\rtimes \mathfrak{f}_n
\end{equation}
есть итерированное разложение в полупрямую сумму, такое что  $\mathfrak{f}_1,\ldots,\mathfrak{f}_n$ одномерны. Пусть также $\varphi\!:\!U(\mathfrak{g})\to K$ --- гомоморфизм $\ptn$-алгебр и выполнены следующие условия:

\emph{(1)}~$K_1,\ldots,K_n$ --- $\ptn$-алгебры Хопфа и имеет место итерированное разложение
\begin{equation}\label{expfsmp}
K\cong(\cdots (K_1 \mathop{\widehat{\#}}K_2)
\mathop{\widehat{\#}}\cdots)\mathop{\widehat{\#}}K_n,
\end{equation}

\emph{(2)}~$\varphi$ согласован с разложением
$$
U(\mathfrak{g})\cong(\cdots (U(\mathfrak{f}_1) \mathop{\#}U(\mathfrak{f}_2))\mathop{\#} \cdots)\mathop{\#} U(\mathfrak{f}_n),
$$
порождённым~\eqref{fgdec}, в смысле, указанном выше. 

Тогда если каждый из гомоморфизмов $U(\mathfrak{f}_i)\to K_i$  имеeт плотный образ и является гомологическим эпиморфизмом, то таков же и $\varphi$.
\end{thm}

В свою очередь, для доказательства теоремы~\ref{AMede2} нужен следующий результат. 

\begin{thm}\label{Assmpr}
\cite[Theorem 4.5]{Ar_smhe}
Пусть $K$ есть  $\mathbin{\widehat{\otimes}}$-алгебра Хопфа, $S$ есть $K$-$\mathbin{\widehat{\otimes}}$-модульная алгебра, и $R$ есть $\mathbb{C}[z]$-модульная алгебра счётной линейной размерности. Предположим, что $\beta\!:\mathbb{C}[z]\to K$ есть гомоморфизм $\mathbin{\widehat{\otimes}}$-алгебр Хопфа и
$\alpha\!:R\to S$ есть гомоморфизм алгебр, являющийся также морфизмом $\mathbb{C}[z]$-модулей. Если оба гомоморфизма $\alpha$ и $\beta$  имеют плотные образы и являются гомологическими эпиморфизмами, то таков же и
$$
\alpha\mathop{\widehat{\#}}\beta\!:R\Smash \mathbb{C}[z]\to S\Psmash K.
$$
\end{thm}

Отметим, что как и в доказательстве теоремы~\ref{Assmpr0} здесь используется то, что $R\Smash \mathbb{C}[z]$ относительно квазисвободна  и имеет $(f_2)$-конечный тип над $R$, см. \cite[Proposition 7.9]{Pir_qfree} и \cite[Proposition 4.10]{Ar_smhe}.

Теперь мы докажем теорему~\ref{AMede2}, а затем теорему~\ref{solvmain}.

\begin{proof}[Доказательство теоремы~\ref{AMede2}]
Рассужения такие же как в \cite[Theorem~4.5]{Ar_smhe}. Тем не менее, здесь приведено полное доказательство.

Рассуждая по индукции, мы покажем, что гомоморфизм
\begin{equation}\label{truAMe}
(\cdots (U(\mathfrak{f}_1) \mathop{\#}U(\mathfrak{f}_2))\mathop{\#} \cdots)\mathop{\#} U(\mathfrak{f}_k)\to (\cdots (K_1 \mathop{\widehat{\#}}K_2)
\mathop{\widehat{\#}}\cdots)\mathop{\widehat{\#}}K_k,
\end{equation}
полученный на каждом шаге итерации при $k=1,\ldots,n$, является гомологическим эпиморфизмом и  имеет плотный образ.

При $k=1$ это утверждение следует из условия. Теперь предположим, что утверждение выполнено при $k-1$. Заметим, что $U(\mathfrak{f}_k)\cong \mathbb{C}[z]$. Так как $\varphi$ согласован с итерированными разложениями, гомоморфизм в \eqref{truAMe} можно записать как
$$
\alpha\mathop{\widehat{\#}}\beta\!:R_{k-1}\Smash \mathbb{C}[z]\to S_{k-1}\Psmash K_k,
$$
где $R_{k-1}$ и $S_{k-1}$ обозначают алгебры, полученные на $(k-1)$-м шаге.  Здесь $\beta\!:\mathbb{C}[z]\to K_k$ --- гомоморфизм  $\mathbin{\widehat{\otimes}}$-алгебр Хопфа, а $\alpha\!:R_{k-1}\to S_{k-1}$ --- это гомоморфизм $\mathbin{\widehat{\otimes}}$-алгебр, который является морфизмом $\mathbb{C}[z]$-модулей.

Оба эпиморфизма $\alpha$ и $\beta$ являются гомологическими и имеют плотные образы (первый по предположению индукции, второй по условию теоремы), а   размерность $R_{k-1}$  счётна.  Таким образом, все предположения теоремы~\ref{Assmpr} выполнены и, следовательно, $R_{k-1}\Smash \mathbb{C}[z]\to S_{k-1}\Psmash K_k$ также является гомологическим эпиморфизмом. Индукция завершена.

Полагая $k=n$, мы заключаем, что $\varphi$ также имеeт плотный образ и является гомологическим эпиморфизмом.
\end{proof}

\begin{proof}[Доказательство теоремы~\ref{solvmain}]
Обозначим  через $\theta$ гомоморфизм пополнения ${\mathscr A}(G)\to{\mathscr A}_{\omega^\infty}(G)$, а через $\tau$ естественным образом определённый гомоморфизм $U(\mathfrak{g})\to {\mathscr A}(G)$ (см.~\eqref{taudef}), и рассмотрим их композицию  $\theta\tau\!:U(\mathfrak{g})\to {\mathscr A}_{\omega^\infty}(G)$.

Из сделанных предположений  о разложениях группы и веса  следует, что ${\mathscr A}_{\omega^\infty}(G)$ может быть представлена в виде итерированного смэш-произведения
$$
(\cdots ({\mathscr A}_{\omega_1^\infty}(F_1) \mathop{\widehat{\#}} {\mathscr A}_{\omega_2^\infty}(F_2))\mathop{\widehat{\#}}\cdots)
\mathop{\widehat{\#}} {\mathscr A}_{\omega_{n-1}^\infty}(F_{n-1})),
$$
см. \cite[теоремы 3.8 и 6.1]{Ar_smash}. Это разложение может записано в виде~\eqref{expfsmp}, где $K_i= \mathfrak{A}_{s_i}$ для некоторого $i$, см. \cite[теорема 6.3]{Ar_smash}, при этом оно согласовано с разложением $\mathfrak{g}$ в итерированную полупрямую сумму~\eqref{fgdec} и соответственным разложением  $U(\mathfrak{g})$ в итерированное смэш-произведение. Так как в силу предложения~\ref{AsHomep} каждый из эпиморфизмов $U(\mathfrak{f}_i)\to\mathfrak{A}_{s_i}$ гомологический, то из теоремы~\ref{AMede2} следует, что $\theta\tau$  также гомологический эпиморфизм.

С другой стороны, так как $G$ односвязна и разрешима, то $\tau\!:U(\mathfrak{g})\to {\mathscr A}(G)$ --- гомологический эпиморфизм \cite[Theorem~8.3]{Pir_stbflat}, а значит, в силу \cite[Proposition~1.8]{T2}, таков же  и $\theta$. Это означает, что $\omega$ --- локализующий вес.
\end{proof}

\subsection*{Завершение доказательства теоремы~\ref{C4maingen}}

Теперь мы можем закончить рассуждения для веса, допускающего разложение, в общем случае.

\begin{proof}[Доказательство теоремы~\ref{C4maingen}]
Обозначим группы  $(\cdots (F_1 \rtimes F_2)\rtimes\cdots)\rtimes F_{n-1}$ и $F_n$  из~\eqref{itdirp0} через $B$ и $L$. Тогда  $G=B\rtimes L$, а $\widetilde\omega_{n-1}$ и $\omega_{n}$ ---  субмультипликативные веса на $B$ и $L$, соответственно. Будучи ограничениями асимптотически симметричного веса, они сами являются асимптотически симметричными. Согласно условию теоремы
\begin{equation*}
\omega(bl)\simeq  \widetilde\omega_{n-1}(b)\omega_n(l)\quad\text{на
$B\times L$} \quad(\text{$b\in B$ и $l\in L$}).
\end{equation*}
В силу \cite[теорема 3.8]{Ar_smash} имеет место изоморфизм $\ptn$-алгебр
$$
{\mathscr A}_{\widetilde\omega_{n-1}^\infty}(B)\mathop{\widehat{\#}} {\mathscr A}_{\omega_n^\infty}(L)\to
{\mathscr A}_{\omega^\infty}(G).
$$
Обозначим ${\mathscr A}(B)\to {\mathscr A}_{\widetilde\omega_{n-1}^\infty}(B)$ и ${\mathscr A}(L)\to {\mathscr A}_{\omega_n^\infty}(L)$ через $\alpha$ и $\beta$ соответственно. Осталось показать, что $\alpha\mathop{\widehat{\#}}\beta$ является гомологическим эпиморфизмом.

Мы воспользуемся теоремой~\ref{Assmpr0}. Действительно, $\beta$ является гомологическим эпиморфизмом в силу предложения~\ref{coredhe}. С другой стороны,  ${\mathscr A}(B)\to {\mathscr A}_{\widetilde\omega_{n-1}^\infty}(B)$ есть гомологический эпиморфизм согласно теореме~\ref{solvmain}. Более того, из теоремы~\ref{redrelhotr}  следует, что ${\mathscr A}(G)$ является относительно гомологически тривиальной над ${\mathscr A}(B)$.  Так как ${\mathscr A}(G)$ изоморфна ${\mathscr A}(B)\mathop{\widehat{\#}} {\mathscr A}(L)$ мы можем применить теорему~\ref{Assmpr0}, из которой получаем, что $\alpha\mathop{\widehat{\#}}\beta$ является гомологическим эпиморфизмом.
\end{proof}

\section{Доказательство теоремы о максимальном весе}
\label{s:proofmainexdi}

Как отмечено выше, доказательство теоремы~\ref{C4mainexdi} включает три шага:
\begin{itemize}
  \item редукция к случаю группы Штейна,
  \item редукция к случаю линейной группы,
  \item доказательство в случае линейной группы.
\end{itemize}

\subsection*{Редукция к случаю группы Штейна}

Сначала мы покажем, что  в теореме~\ref{C4mainexdi} достаточно рассматривать только (связные) группы Штейна.

Согласно Моримото \cite[Theorem~1]{Mo65} всякая связная комплексная группа Ли~$G$ содержит наименьшую замкнутую нормальную подгруппу $\mathrm{M}$, такую что $G/\mathrm{M}$ является группой Штейна. Будем называть её \emph{подгруппой Моримото}.

\begin{pr}\label{Steiniz}
Гомоморфизм $\pi\!:{\mathscr A}(G)\to{\mathscr A}(G/\mathrm{M})$, индуцированный факторотображением $G\to G/M$, является топологическим изоморфизмом.
\end{pr}

Нам понадобится следующая лемма, которая проверяется непосредственно.

\begin{lm}\label{consthol}
Пусть $N$ --- замкнутая  нормальная комплексная подгруппа в комплексной группе Ли~$G$.
Линейное отображение $\mathcal{O}(G/N)\to\mathcal{O}(G)$, индуцированное факторотображением $G\to G/N$,  топологически инъективно, а его образ совпадает с пространством голоморфных функций, постоянных на каждом смежном классе~$N$.
\end{lm}

\begin{proof}[Доказательство предложения~\ref{Steiniz}]
Из \cite[Theorem~1]{Mo65} следует, что всякая голоморфная функция на $G$ постоянна на смежных классах подгруппы Моримото~$\mathrm{M}$. Из леммы~\ref{consthol} получаем, что линейное отображение $\pi'\!:\mathcal{O}(G/\mathrm{M})\to\mathcal{O}(G)$, двойственное к $\pi$, с одной стороны, сюръективно, а с другой, топологически инъективно, т.е. является топологическим изоморфизмом. Тем самым $\pi$ --- также топологический изоморфизм.
\end{proof}

\subsection*{Редукция к случаю линейной группы}

Теперь мы покажем, что можно предполагать, что группа линейна. Следующее утверждение --- основной результат работы~\cite{Ar_lin} (собственно говоря, упомянутая статья и была написана автором, чтобы в дальнейшем использовать свойства линеaризатора при исследовании гомологических эпиморфизмов).

\begin{thm}\label{linarizc}
\cite[теорема 1]{Ar_lin} 
Пусть $G$ --- связная комплексная группа Ли, $\Lambda$ --- её линеaризатор, а $\mathrm{M}$ --- её подгруппа  Моримото. Тогда $\Lambda$ изоморфен $\mathrm{M}\times(\mathbb{C}^\times)^k$ для некоторого $k\in\mathbb{Z}_+$. 
\end{thm}

\begin{co}\label{linarizcco}
Пусть $\Lambda$ --- линеaризатор связной комплексной группы Ли. 
Тогда ${\mathscr A}(\Lambda)$ гомологически тривиальна.
\end{co}

\begin{proof}
Из теоремы~\ref{linarizc} и предложения~\ref{Steiniz} следует, что ${\mathscr A}(\Lambda)$ изоморфна ${\mathscr A}((\mathbb{C}^\times)^k)$ для некоторого $k\in\mathbb{Z}_+$. Будучи комплексификацией  тора, $(\mathbb{C}^\times)^k$  линейно комплексно редуктивна. Тем самым ${\mathscr A}((\mathbb{C}^\times)^k)$  гомологически тривиальна в силу теоремы~\ref{Aredhtr}.
\end{proof}

Следствие~\ref{linarizcco} будет использовано в доказательстве предложения~\ref{qLaloca}. 
А теперь мы покажем, что при изучении некоторых гомологических вопросов можно игнорировать всякую замкнутую  нормальную комплексную подгруппу $N$ такую, что ${\mathscr A}(N)$ гомологически тривиальна.

Если $N$ --- замкнутая  нормальная комплексная подгруппа в комплексной группе Ли~$G$, то
обозначим через ${\mathscr A}(N)_0$ ядро коединицы ${\mathscr A}(N)$ (т.е. множество всех аналитических функционалов, равных~$0$ на постоянных функциях).
Так как гомоморфизм ${\mathscr A}(N)\to {\mathscr A}(G)$, порождённый вложением $N\to G$, топологически инъективен (см. теорему~\ref{subgtopin} в приложении), то можно рассматривать ${\mathscr A}(N)$ и ${\mathscr A}(N)_0$ как замкнутые подмножества в ${\mathscr A}(G)$. По сути дела, эта информация немного избыточна, так как в лемме~\ref{kerpiovl} и предложении~\ref{AGNAGpr} ниже, мы используем лишь то, что гомоморфизм ${\mathscr A}(N)\to {\mathscr A}(G)$ снабжает ${\mathscr A}(G)$  структурой ${\mathscr A}(N)$-$\ptn$-модуля. Однако поскольку топологическая инъективность имеет место, мы не будем её игнорировать.

\begin{lm}\label{kerpiovl}
Пусть $N$ --- замкнутая  нормальная комплексная подгруппа в комплексной группе Ли~$G$ и пусть
 $\pi\!:{\mathscr A}(G)\to{\mathscr A}(G/N)$  ---  гомоморфизм  $\ptn$-алгебр, индуцированный
факторотображением $G\to G/N$. Тогда
$$
\Ker\pi =\,\overline{{\mathscr A}(G)\,{\mathscr A}(N)_0}\,=\,\overline{{\mathscr A}(N)_0\, {\mathscr A}(G)},
$$
где черта означает замыкание.
\end{lm}
\begin{proof}
Мы используем тот же трюк, что и доказательстве \cite[теорема~2.1]{Ar_smash} --- редукцию к групповым алгебрам (которые плотны в алгебрах аналитических функционалов в силу леммы~\ref{deltden}).

Обозначим очевидно определённый гомоморфизм $\mathbb{C}[G]\to\mathbb{C}[G/N]$ через $\pi_d$ и покажем, что образ $\Ker \pi_d$ при вложении $\mathbb{C}[G]\to {\mathscr A}(G)$ плотен в $\Ker \pi$. Для этого достаточно проверить, что двойственное отображение  $(\Ker \pi)' \to (\Ker \pi_d)'$ инъективно.
Действительно, в силу леммы~\ref{consthol} $(\Ker \pi)' $ можно отождествить с $\mathcal{O}(G)/\Image \pi'$ и аналогично $(\Ker \pi_d)'$ можно отождествить с $\mathbb{C}^G/\Image \pi_d'$. Если функция из $\mathcal{O}(G)$ содержится в $\Image \pi_d'$, то она очевидно содержится в $\Image \pi'$. Этим доказана инъективность $(\Ker \pi)' \to (\Ker \pi_d)'$.

Обозначим ядро коединицы групповой алгебры $\mathbb{C}[N]$ через $\mathbb{C}[N]_0$. Оно состоит из тех линейных комбинаций дельта-функций, суммы коэффициентов которых равны~$0$.
Проверим, что $\Ker \pi_d=\mathbb{C}[G]\,\mathbb{C}[N]_0$. Действительно, если $\sum_{h\in G} c_h\delta_h\in \Ker\pi_d$, то для всякого $g\in G$ имеем $\sum_{h\in gN}c_h=0$ и тем самым $\sum_{h\in gN}c_h\delta_{g^{-1}h}\in\mathbb{C}[N]_0$. Отсюда следует, что $\Ker \pi_d\subset\mathbb{C}[G]\,\mathbb{C}[N]_0$. Обратное включение очевидно.

В силу леммы~\ref{deltden} образы  $\mathbb{C}[G]$ и $\mathbb{C}[N]_0$ плотны в ${\mathscr A}(G)$ и ${\mathscr A}(N)_0$ соответственно. Используя равенство $\Ker \pi_d=\mathbb{C}[G]\,\mathbb{C}[N]_0$ и плотность образа $\Ker \pi_d$ в $\Ker \pi$, заключаем отсюда, что $\Ker\pi =\,\overline{{\mathscr A}(G)\,{\mathscr A}(N)_0}$. Вторая формула получается аналогично.
\end{proof}

Нам понадобится следующий хорошо известный изоморфизм, см. доказательство в нужной нам общности в
\cite[Proposition~3.1]{Pi09}.
\begin{pr}\label{wmfor}
Пусть $A$ есть $\ptn$-алгебра, а $I$ есть её замкнутый левый идеал. Тогда
для каждого $X\in \rmod A$ существует топологический изоморфизм
$$
X\ptens{A}(A/I)^\sim \cong(X/\,\overline{X\cdot I})^\sim
$$
(здесь тильда обозначает пополнение),
определённый формулой $$x \otimes (a + I)\mapsto x\cdot a + X \cdot I.$$
\end{pr}

Рассмотрим на ${\mathscr A}(G/N)$ структуру ${\mathscr A}(G)$-$\ptn$-бимодуля, порождённую гомоморфизмом $\pi\!:{\mathscr A}(G)\to{\mathscr A}(G/N)$.

\begin{pr}\label{AGNAGpr}
Пусть $N$ --- замкнутая  нормальная комплексная подгруппа в комплексной группе Ли~$G$. Если
${\mathscr A}(N)$ гомологически тривиальна, то ${\mathscr A}(G/N)$ является проективным и как левый и как правый
${\mathscr A}(G)$-$\ptn$-модуль.
\end{pr}
\begin{proof}
Заметим, что $\pi\!:{\mathscr A}(G)\to{\mathscr A}(G/N)$   сюръективен и открыт
\cite[Proposition~A.7]{AHHFG}.
Более того, согласно лемме~\ref{kerpiovl} имеем $\Ker \pi
=\overline{{\mathscr A}(G)\,{\mathscr A}(N)_0}$.  Тем самым
$$
{\mathscr A}(G/N)\cong {\mathscr A}(G)/\,\overline{{\mathscr A}(G)\,{\mathscr A}(N)_0},
$$
и пространство в правой части полно. (Мы рассматриваем  ${\mathscr A}(N)_0$ и ${\mathscr A}(N)$ как замкнутые подалгебры в ${\mathscr A}(G)$, см. теорему~\ref{subgtopin} в приложении.) Будучи одномерным, пространство ${\mathscr A}(N)/{\mathscr A}(N)_0$ также полно. Полагая  $A={\mathscr A}(N)$, $I={\mathscr A}(N)_0$ и $X={\mathscr A}(G)$ в предложении~\ref{wmfor}, получаем
$$
{\mathscr A}(G)/\,\overline{{\mathscr A}(G)\,{\mathscr A}(N)_0}\cong
{\mathscr A}(G)\ptens{{\mathscr A}(N)}({\mathscr A}(N)/{\mathscr A}(N)_0),
$$
так как нет необходимости в пополнении указанных выше двух пространств.
Объединяя обе формулы с изоморфизмом
${\mathscr A}(N)/{\mathscr A}(N)_0\cong\mathbb{C}$, заключаем, что
${\mathscr A}(G/N)\cong {\mathscr A}(G)\ptens{{\mathscr A}(N)}\mathbb{C}$.

Так как ${\mathscr A}(N)$ гомологически тривиальна, то $\mathbb{C}$ проективен в ${\mathscr A}(N)\lmod$.
Остаётся воспользоваться тем очевидным фактом, что  функтор расширения скаляров ${\mathscr A}(G)\ptens{{\mathscr A}(N)}(-)$ переводит
свободные объекты в свободные, а значит, проективные в проективные.

Доказательство для правых модулей аналогично.
\end{proof}

Теперь мы можем доказать утверждение о редукции теоремы~\ref{C4mainexdi} к случаю связной линейной группы.

\begin{pr}\label{qLaloca}
Пусть~$G$ --- связная комплексная группа Ли, $\Lambda$ --- её линеaризатор, а $\omega$ --- локализующий вес на~$G/\Lambda$. Тогда эпиморфизм ${\mathscr A}(G)\to{\mathscr A}_{\omega^\infty}(G/\Lambda)$ является гомологическим.
\end{pr}
\begin{proof}
Будучи сюръективным \cite[Proposition A.7]{AHHFG}, гомоморфизм ${\mathscr A}(G)\to {\mathscr A}(G/\Lambda)$ является эпиморфизмом. 

Согласно следствию~\ref{linarizcco} алгебра ${\mathscr A}(\Lambda)$ гомологически тривиальна. Поэтому мы можем воспользоваться предложением~\ref{AGNAGpr}, из которого следует, что 
${\mathscr A}(G/\Lambda)$ есть проективный левый ${\mathscr A}(G)$-$\ptn$-модуль. Тогда непосредственно из определения получаем, что эпиморфизм ${\mathscr A}(G)\to {\mathscr A}(G/\Lambda)$ является гомологическим. Так как $\omega$ --- локализующий вес, то ${\mathscr A}(G/\Lambda)\to{\mathscr A}_{\omega^\infty}(G/\Lambda)$ также является гомологическим эпиморфизмом. Осталось
воспользоваться тем, что композиция двух гомологических эпиморфизмов есть гомологический эпиморфизм.
\end{proof}

\subsection*{Доказательство теоремы~\ref{C4mainexdi} в общем случае}
Теперь мы можем завершить доказательство ключевой теоремы.

\begin{proof}[Доказательство теоремы~\ref{C4mainexdi}]
(1)~Сначала предположим, что $G$ линейна. Согласно \cite[теорема~4.4]{Ar_smash} при указанных условиях группа $G$ и вес $\omega_{\max}$ допускают разложения \eqref{itdirp0} и~\eqref{weidec} соответственно, причём веса $\omega_i$ имеют вид либо $z\mapsto 1+|z|$ либо $z\mapsto \exp(|z|^{1/s})$ для некоторого $s\ge 1$. В силу леммы~\ref{1dimdsp} и предложения~\ref{AsHomep} все эти веса являются локализующими. Поэтому из теоремы~\ref{C4maingen} следует, что  $\omega_{\max}$ --- локализующий.

(2)~Пусть~$G$ --- произвольная связная  комплексная группа Ли, а $\Lambda$ --- её линеaризатор. В силу части~(1) заключаем, что $\omega_{\max}$ --- локализующий вес на~$G/\Lambda$. Тогда согласно предложению~\ref{qLaloca} эпиморфизм ${\mathscr A}(G)\to{\mathscr A}_{\omega^\infty}(G/\Lambda)$ является гомологическим.
\end{proof}

Напомним, что основная теорема~\ref{C4main} является простым следствием теоремы~\ref{C4mainexdi},
см. \S\,\ref{s:mainfirst}.

Как видно из предыдущих результатов, локализующие веса существуют в изобилии. Рассмотрим также пример веса, который не является локализующим,

\begin{exm}\label{nonlovw}
Рассмотрим на $G=\mathbb{C}$ вес $\omega(z)=1+\log|z|$. Из \cite[лемма 2.10]{ArAMN} следует, что в этом случае ${\mathscr A}_{\omega^\infty}(G)$ --- одномерная алгебра. Применяя  ${\mathscr A}_{\omega^\infty}(G)\ptens{{\mathscr A}(G)}(-)$ к резольвенте
\begin{equation*}
0 \leftarrow \mathbb{C} \xleftarrow{\varepsilon} {\mathscr A}(\mathbb{C})\xleftarrow{\delta_0}{\mathscr A}(\mathbb{C})\leftarrow 0.
\end{equation*}
(см. доказательство предложения~4.4 из \cite{Ar_smhe}) получаем последовательность
$$
0 \longleftarrow \mathbb{C} \longleftarrow \mathbb{C} \longleftarrow \mathbb{C} \leftarrow 0,
$$
которая не точна. Тем самым эпиморфизм ${\mathscr A}(G)\to{\mathscr A}_{\omega^\infty}(G)$  не является гомологическим, а значит $\omega$ --- локализующим.
\end{exm}

\section{Доказательство теоремы о пополнении универсальной обертывающей}
\label{s:proofUg}

Доказательство теоремы~\ref{AMheifnew} (о пополнении универсальной обертывающей алгебры), приведённое ниже, существенно опирается на теорему~\ref{AMede2} и на результаты Пирковского из \cite{Pi4,Pir_stbflat}.

\begin{proof}[Доказательство теоремы~\ref{AMheifnew}]
(1)~Необходимость. Предположим, что эпиморфизм $U(\mathfrak{g})\to \widehat U(\mathfrak{g})^{\mathrm{PI}}$  является гомологическим и покажем, что $\mathfrak{g}$ разрешима.

Рассуждение аналогично доказательству \cite[Theorem~3.6]{Pi4}, где
используется разложение $\widehat U(\mathfrak{g})$ в аналитическое смэш-произведение, соответствующее разложению Леви алгебры $\mathfrak{g}$. Алгебра $\widehat U(\mathfrak{g})^{\mathrm{PI}}$  также может быть представлена как аналитическое смэш-произведение. Действительно, заметим сначала, что в силу предложения~\ref{evUgAG} $\widehat U(\mathfrak{g})^{\mathrm{PI}}\cong \widehat{\mathscr A}(G)^{\mathrm{PI}}$, где $G$ --- соответствующая односвязная группа Ли.

Пусть $\mathfrak{g}=\mathfrak{r}\rtimes \mathfrak{s}$ --- разложение Леви ($\mathfrak{r}$ --- радикал, а $\mathfrak{s}$ полупроста). Ему соответствует полупрямое произведение $G=R\rtimes S$, где $G$, $R$ и $S$ --- односвязные  комплексные группы Ли, для которых соответственно  $\mathfrak{g}$, $\mathfrak{r}$ и $\mathfrak{s}$  есть их алгебры Ли \cite[глава~6, \S\,2, с.\,197, теорема 2]{VOsem}. Тогда $R$ разрешима и односвязна, а $S$ линейно комплексно редуктивна, будучи полупростой. В частности, это означает, что $G$ линейна и мы можем применить результаты о разложении в аналитическое смэш-произведение. Из \cite[теорема 6.6]{Ar_smash} и существования изоморфизма $\widehat{\mathscr A}(S)\cong \widehat U(\mathfrak{s})$ следует, что
$$
\widehat U(\mathfrak{g})^{\mathrm{PI}}\cong A \mathop{\widehat{\#}} \widehat U(\mathfrak{s}),
$$
где 
\begin{equation}\label{solbPIdec}
A=(\cdots (\mathbb{C}[[x_1]] \mathop{\widehat{\#}} \cdots \mathbb{C}[[x_p]])\mathop{\widehat{\#}} \mathcal{O}(\mathbb{C}))\mathop{\widehat{\#}}\cdots\mathcal{O}(\mathbb{C}).
\end{equation}
Так же, как в доказательстве леммы~3.4 из \cite{Pi4}, получаем, что  ${\mathop{\mathrm{Tor}}\nolimits}_n^{\widehat U(\mathfrak{g})^{\mathrm{PI}}}(\mathbb{C},A)=0$  для всех $n>0$.

С другой стороны,  рассуждая так же, как и в доказательстве леммы~3.5 из \cite{Pi4} с заменой алгебры, которая там обозначена через $\widetilde U(\mathfrak{r})$ на $A$, получаем, что ${\mathop{\mathrm{Tor}}\nolimits}_k^{U(\mathfrak{g})}(\mathbb{C},A)\neq 0$, где $k$ --- размерность $\mathfrak{s}$. (Всё, что требуется от $A$ --- тот факт, что $\mathfrak{s}$ действует на ней непрерывными дифференцированиями и наличие канонической аугментации.)  Отсюда вытекает, что $k=0$, т.е. $\mathfrak{s}$ тривиальна, а значит $\mathfrak{g}$ разрешима.

(2)~Достаточность. Предположим, что $\mathfrak{g}$ разрешима. Тогда 
из \cite[теорема 6.6]{Ar_smash} следует, что существует итерированное разложение~\eqref{solbPIdec} алгебры $A=\widehat U(\mathfrak{g})^{\mathrm{PI}}$,
согласованное с некоторым итерированным разложением в полупрямую сумму
\begin{equation*}
\mathfrak{g}=((\cdots (\mathfrak{f}_1 \rtimes \mathfrak{f}_2)\rtimes\cdots)\rtimes \mathfrak{f}_n
\end{equation*}
и таким, что  $\mathfrak{f}_1,\ldots,\mathfrak{f}_n$ одномерны. Поскольку $U(\mathfrak{f}_i)\cong\mathbb{C}[x]$, а каждый из эпиморфизмов $\mathbb{C}[x]\to \mathbb{C}[[x]]$ и $\mathbb{C}[x]\to\mathcal{O}(\mathbb{C})$ является гомологическим (лемма~\ref{AsHomep}) и имеет плотный образ, из теоремы~\ref{AMede2} следует, что $U(\mathfrak{g})\to \widehat U(\mathfrak{g})^{\mathrm{PI}}$ также является гомологическим эпиморфизмом.
\end{proof}

\appendix

\section{Топологические вложения пространств аналитических функционалов}
\label{a:topem}

Мы ссылаемся на теорему, приведённую в этом приложении, в доказательстве предложения~\ref{AGNAGpr}. Как указано в обсуждении перед леммой~\ref{kerpiovl}, она не является необходимой частью рассуждений, но, тем не менее, облегчает понимание поведения алгебр аналитических функционалов на подгруппах. Сама по себе теорема представляет самостоятельный интерес и доказана здесь в более общем контексте пространств Штейна.

Напомним, что всякий морфизм комплексных аналитических пространств $$(S_1,\mathcal{O}_{S_1})\to (S_2,\mathcal{O}_{S_2})$$  индуцирует гомоморфизм алгебр глобальных сечений $\rho\!:\mathcal{O}(S_2)\to \mathcal{O}(S_1)$, а всякая алгебра сечений снабжена канонической топологией, превращающей её в алгебру Фреше. Более того, $\rho$ непрерывен относительно канонических топологий. Подробности см., например, в начале \S\,1 из~\cite{Fo67}. Отметим только, что топологии на алгебрах сечений сначала определяются  для достаточно малых открытых подмножеств, а затем переносятся на произвольные с помощью аксиомы склейки. Непрерывность также сначала доказывается для достаточно малых открытых подмножеств, а затем распространяется на общий случай с использованием универсального свойства уравнителя, фигурирующего в аксиоме склейки.

Так же как и в случае многообразия, положим ${\mathscr A}(S_1)\!:=\mathcal{O}(S_1)'$ для комплексного аналитического пространства $(S_1,\mathcal{O}_{S_1})$ (здесь штрих обозначает сильное двойственное пространство).
Нас интересуют свойства двойственного линейного отображения $\rho'\!:{\mathscr A}(S_1)\to {\mathscr A}(S_2)$ в случае, когда $(S_1,\mathcal{O}_{S_1})\to (S_2,\mathcal{O}_{S_2})$  является замкнутым вложением пространств Штейна (см. определение и основные свойства в \cite[Chapter~1, \S\,2, п.\,2, c.\,15, и п.\,7, c.\,20--21]{GR2}).

\begin{thm}\label{subgtopin}
Пусть $(S_1,\mathcal{O}_{S_1})\to (S_2,\mathcal{O}_{S_2})$  есть замкнутое вложение пространств Штейна. Тогда линейное отображение $\rho'\!:{\mathscr A}(S_1)\to {\mathscr A}(S_2)$, им порождённое, топологически инъективно.
\end{thm}

Следующая лемма следует непосредственно из \cite[Proposition 26.23]{MV}.
\begin{lm}\label{topxyprxpr}
Пусть $X$ --- пространство Фреше, а $Y$ --- его замкнутое подпространство, причём $X/Y$ является монтелевским. Тогда  отображение  $\sigma'\!:(X/Y)'\to X'$ сильных двойственных пространств, двойственное к проекции $\sigma\!:X\to X/Y$,  является топологически инъективным.
\end{lm}

\begin{proof}[Доказательство теоремы~\ref{subgtopin}]
Можно отождествить $S_1$ c замкнутым аналитическим подпространством в $S_2$. В частности, это значит, что $\mathcal{O}_{S_1}$  изоморфен факторпучку пучка $\mathcal{O}_{S_2}$. Так как $(S_2,\mathcal{O}_{S_2})$ --- пространство Штейна, мы можем применить теорему~B Картана и заключить, что гомоморфизм $\rho\!:\mathcal{O}(S_2)\to \mathcal{O}(S_1)$ сюръективен.

Далее, $\mathcal{O}(S_1)$  является ядерным пространством Фреше, см., например, \cite[Lemma A.4]{AHHFG}, а всякое ядерное пространство Фреше является монтелевским, см., например, \cite[глава III, \S\,7.2, с.\,131, следствие~2]{Schae}. 

Так как $\mathcal{O}(S_2)$ и $\mathcal{O}(S_1)$ --- пространства Фреше, то из теоремы об открытом отображении следует, что $\mathcal{O}(S_1)\cong \mathcal{O}(S_2)/\Ker \rho$. Тем самым мы можем отождествить $\rho'\!:{\mathscr A}(S_1)\to {\mathscr A}(S_2)$ с $(\mathcal{O}(S_2)/\Ker \rho)'\to \mathcal{O}(S_2)'$. Итак,  положив $X=\mathcal{O}(S_1)$ и $Y=\Ker \rho$ в лемме~\ref{topxyprxpr}, получаем, что отображение $\rho'$  топологически инъективно.
\end{proof}


\begin{thebibliography}{99}

\bibitem{Ak08}
С.\,С.~Акбаров,  \emph{Голоморфные функции экспоненциального типа и двойственность для групп Штейна с алгебраической связной компонентой единицы}, Фундамент. и прикл. матем., 14:1 (2008), 3--178; J. Math. Sci., 162:4 (2009), 459--586.

\bibitem{ArAMN}
O.\,Yu.~Aristov, \emph{Arens--Michael envelopes of nilpotent Lie
algebras, functions of exponential type, and homological
epimorphisms}, Тр. ММО, 81, № 1, МЦНМО, М., 2020, 117--136, Trans. Moscow Math. Soc. (2020), 97--114.

\bibitem{ArAnF}
O.\,Yu.~Aristov, \emph{Holomorphic functions of exponential type on
connected complex Lie groups}, J. Lie Theory 29:4 (2019),
1045--1070, arXiv:1903.08080.

 \bibitem{ArHR}
 O.\,Yu.~Aristov, \emph{On holomorphic reflexivity conditions for complex Lie groups}, Proc. Edinb. Math. Soc.,  64:4 (2021), 800--821.

\bibitem{ArRC}
О.\,Ю.~Аристов, \emph{Соотношение “коммутатор равен функции” в банаховых алгебрах}, Матем. заметки, 109:3 (2021), 323--337; Math. Notes, 109:3 (2021), 323--334.

\bibitem{Ar_lin}
О.\,Ю.~Аристов, \emph{Строение линеаризатора связной комплексной группы Ли}, Сиб. матем. журн., 64:2 (2023), 276--280.

\bibitem{ArPiLie}
O.\,Yu.\,Aristov, \emph{When a completion of the universal enveloping algebra is a Banach PI-algebra?}, Bull. Aust. Math. Soc, 107:3 (2023), 493--501, arXiv: 2204.07393.

\bibitem{ArLfd}
O.\,Yu.\,Aristov,
\emph{Length functions exponentially distorted on subgroups of complex Lie groups}, European J. Math., 9 (2023), 60,  arXiv: 2208.12667.

\bibitem{AHHFG}
O.\,Yu.~Aristov, \emph{Holomorphically finitely generated Hopf algebras and quantum Lie groups}, Moscow Math. J., 24:2 (2024), 145--180,  arXiv:2006.12175.


\bibitem{Ar_smash}
О.\,Ю.~Аристов, \emph{Разложение алгебры аналитических функционалов на связной комплексной группе Ли и её пополнений в итерированные аналитические смэш-произведения}, Алгебра и анализ, 36:4 (2024), 1--37, arXiv: 2209.04192.

\bibitem{Ar_envPG}
O.\,Yu.~Aristov, \emph{Envelopes in the class of Banach algebras of polynomial growth and $C^\infty$-functions of a finite number of free variables}, J. Funct. Anal., 289 (2025), 111117, arXiv: 2401.10199.

\bibitem{Ar_smhe}
O.\,Yu.~Aristov, \emph{The Arens--Michael envelope of a solvable Lie algebra is a homological epimorphism}, Front. Math., DOI: 10.1007/s11464-024-0114-5, arXiv: 2404.19433.

\bibitem{Ar_dim}
О.\,Ю.~Аристов, \emph{Гомологические размерности алгебр аналитических функционалов и их пополнений}, Изв. РАН. Сер. матем. (в печати), arXiv: 2510.26249.

\bibitem{AP}
O.\,Yu.~Aristov, A.\,Yu.~Pirkovskii, \emph{Open embeddings and pseudoflat epimorphisms}, J. Math. Anal. Appl. 485 (2020) 123817.

\bibitem{BBK18}
F.\,Bambozzi, O.\,Ben-Bassat, K.\,Kremnizer, \emph{Stein domains in Banach algebraic geometry}, J.
Funct. Anal. 274:7 (2018),  1865--1927.

\bibitem{BS01}
D.~Belti\c{t}\u{a}, M.~\c{S}abac, \emph{Lie Algebras of Bounded
Operators},   Operator Theory: Advances and Applications, 120.
Birkh\"{a}user Verlag, Basel, 2001.

\bibitem{Bi20}
B.\,Bilich, \emph{Taylor spectrum for modules over Lie algebras}, Funktsional. Anal. i Prilozhen., 56:3 (2022), 3--15; Funct. Anal. Appl. 56:3  (2022), 159--168,  arXiv:2108.12415.

\bibitem{BL93}
E.~Boasso, A.~Larotonda,  \emph{A  spectral  theory  for  solvable  Lie  algebras  of operators},  Pacific  J.  Math.
158 (1993), 15--22.

\bibitem{Bou}
N.\,Bourbaki, \emph{Elements of mathematics. Lie groups and Lie algebras. Part I: Chapters 1--3}, Addison-Wesley/Hermann, Paris, 1975.

\bibitem{Co85}
A.\,Connes, \emph{Non-commutative differential geometry}, Inst. Hautes \'{E}tudes Sci. Publ. Math.
(1985), 41--144.

\bibitem{Co08}
Y.\,de Cornulier, \emph{Dimension of asymptotic cones of Lie
groups}, J. of Topology, 1 (2008), 342--361.

\bibitem{Dix}
J.\,Dixmier, \emph{Enveloping algebras}, North-Holland Publ., 1977.

\bibitem{Do05}
А.\,А.~Досиев, \emph{Когомологии пучков алгебр Фреше и спектральная теория}, Функц. анализ и его прил., 39:3 (2005), 76--80; Funct. Anal. Appl., 39:3 (2005), 225--228.

\bibitem{Do09C}
Dosi A., \emph{Fr\'echet sheaves and Taylor spectrum for supernilpotent Lie algebra of operators}, Mediterr. J. Math. 6 (2009), 181--201.

\bibitem{Do09}
A.\,A.~Dosiev (Dosi), \emph{Local left invertibility for operator tuples and
noncommutative localizations}. J. K-Theory 4:1 (2009), 163--191.

\bibitem{Do10C}
А.\,А.~Доси, \emph{Спектр Тейлора и трансверсальность для операторной алгебры Гейзенберга}, Матем. сб., 201:3 (2010), 39--62; A.\,A.~Dosi, \emph{The Taylor spectrum and transversality for a Heisenberg algebra of operators}, Sb. Math., 201:3 (2010), 355--375.

\bibitem{Do10A}
A.\,A.~Dosi, \emph{Taylor functional calculus for supernilpotent Lie algebra of operators, J. of Operator Th.}, 63:1 (2010), 191--216.

\bibitem{Do10B}
A.\,A.~Dosiev (Dosi), \emph{Formally-radical functions in elements of a nilpotent Lie algebra and noncommutative localizations}, Algebra Colloq., 17, Sp. Iss. 1 (2010), 749--788.

\bibitem{Fa93}
 A.\,S.~Fainstein, \emph{Taylor joint spectrum for  families of operators generating nilpotent Lie algebras},  J.  Operator  Theory 29 (1993), 3--27.

\bibitem{Fo67}
O.~Forster, \emph{Zur Theorie der Steinschen Algebren und Moduln},
Math.~Z. 97 (1967), 376--405.

\bibitem{GR2}
H.~Grauert, R.~Remmert,  \emph{Coherent analytic sheaves}, Springer, Heidelberg 1984.

\bibitem{He81}
А.\,Я.~Хелемский, \emph{Гомологические методы в голоморфном исчислении от нескольких операторов в банаховом пространстве, по Тейлору}, УМН, 36:1(217) (1981), 127--172; Russian Math. Surveys, 36:1 (1981), 139--192.

\bibitem{X1}
A.\,Я.~Хелемский, \emph{Гомология в банаховых и топологических алгебрах}, М., МГУ, 1986.

\bibitem{HiNe}
J.~Hilgert, K.-H.~Neeb, \emph{Structure and geometry of Lie
groups}, Springer, 2011.

\bibitem{Ji92}
R.\,Ji, \emph{Smooth dense subalgebras of reduced group, and $C^*$-algebras, Schwartz cohomology of groups, and cyclic
cohomology}, J. Funct. Anal. 107:1 (1992), 1--33.

\bibitem{KKR16}
A.\,Kanel-Belov, Y.\,Karasik, L.\,H.~Rowen, \emph{Computational aspects of
polynomial identities, Volume~I, Kemer's theorems}, 2nd ed., 2016.

\bibitem{Kot1}
G.\,K\"{o}the, \emph{Topological vector spaces I},  Springer, New York, 1969.

\bibitem{Li70}
Г.\,Л.~Литвинов, \emph{Групповые алгебры аналитических функционалов и их представления}, Докл. АН СССР, 190:4 (1970), 769--771.

\bibitem{MV}
R.\,Meise, D.\,Vogt, \emph{Introduction to functional analysis}, Oxford University Press 1997.


\bibitem{Me04}
R.\,Meyer, \emph{Embeddings of derived categories of bornological modules},
arXiv:041059  (2004).

\bibitem{Me06}
R.~Meyer, \emph{Combable groups have group cohomology of polynomial growth}, The Quarterly Journal of Mathematics, 57:2 (2006), 241--261.

\bibitem{Mo65}
A.~Morimoto, \emph{Non-compact complex Lie groups without non-constant holomorphic
functions}, Proc. of the Conf. on Complex Analysis, Minneapolis 1964, Springer, 1965, 256--272.

\bibitem{VOsem}
Э.\,Б.~Винберг, А.\,Л.~Онищик, \emph{Семинар по группам Ли и алгебраическим группам}, 2 изд. УРСС, 1995.

\bibitem{Os02}
D.\,V.\,Osin, \emph{Exponential radicals of solvable Lie groups}, J.
of Algebra, 248 (2002), 790--805.

\bibitem{Pi99}
A.\,Yu.~Pirkovskii, \emph{On certain homological properties of Stein algebras}, Functional analysis, 3. J. Math.Sci. 95:6 (1999), 2690--2702.

\bibitem{Pi4}
A.\,Yu.~Pirkovskii, \emph{Arens--Michael enveloping algebras and analytic smash
products}, Proc. Amer. Math. Soc. 134 (2006), 2621--2631.

\bibitem{Pir_stbflat}
A.\,Yu.~Pirkovskii, \emph{Stably flat completions of universal enveloping algebras},
Dissertationes Math. (Rozprawy Math.) 441 (2006), 1--60.

\bibitem{Pir_qfree}
A.\,Yu.~Pirkovskii, \emph{Arens--Michael envelopes, homological epimorphisms, and relatively
quasi-free algebras}, Trans. Moscow Math. Soc. 2008, 27--104.

\bibitem{Pi09}
A.\,Yu.~Pirkovskii,
\emph{Flat cyclic Fr\'echet modules, amenable Fr\'echet algebras, and approximate identities}, Homology, Homotopy and Applications, 11:1, 2009, 81--114.

\bibitem{Pir_co1}
A.\,Yu.~Pirkovskii, \emph{A note on relative homological epimorphisms of topological algebras},
arXiv:2104.13716.

\bibitem{Pir_co2}
А.\,Ю.~Пирковский, \emph{Письмо в редакцию}, Тр. ММО, 82:2 (2021),  393--394.

\bibitem{Sa02}
Y.\,Saburi, \emph{On a generalization of the Laurent expansion}, in Microiocal Analysis and
Complex Fourier Analysis, eds. T.\,Kawai, K.\,Fujita, World Scientific, 2002.

\bibitem{Schae}
Х.\,Шефер, \emph{Топологические векторные пространства}, М.: Мир, 1971.

\bibitem{Si55}
Ж.\,Себаштьян-и-Силва, \emph{О некоторых классах локально выпуклых пространств, важных в приложениях}, Математика, 1:1 (1957), 60--77; Rendiconti di matematica e delle sue applicazioni, Roma, 14:5 (1955), 388--410.

\bibitem{T70}
J.\,L.~Taylor, \emph{A joint spectrum for several commuting operators}, J. Functional Analysis
6 (1970), 172--191.

\bibitem{T1}
J.\,L.~Taylor, \emph{Homology and cohomology for topological algebras}, Adv. Math.
\textbf{9} (1972), 137--182.

\bibitem{T2}
J.\,L.~Taylor, \emph{A general framework for a multi-operator functional calculus},
Adv. Math. 9 (1972), 183--252.

\bibitem{Ta73}
J.\,L.~Taylor, \emph{Functions of several noncommuting variables}, Bull. Amer. Math. Soc. 79 (1973), 1--34.

\bibitem{Yo57}
K.\,Yoshinaga, \emph{On a locally convex space introduced by J.S.E. Silva}, Journal of Science
of the Hiroshima University. Series A (Mathematics, Physics, Chemistry) 21:2 (1957), 89--98.
\end{thebibliography}
\end{document}